\documentclass[a4paper, 10pt]{article}

\usepackage{amsfonts, amsthm}
\usepackage{graphics}
\usepackage{caption2}
\usepackage{amssymb}
\usepackage{psfrag}
\usepackage{amsmath}
\usepackage{graphicx}
\usepackage{multirow}
\usepackage{subfigure}

\usepackage{algorithm,algpseudocode,float}



\def\p{\partial}
\def\DD{\displaystyle}
\def\R{\mathbb{R}}
\def\ii{\mathbf{i}}

\def\l{\,l}
\def\bx{{\bf x}}
\def\by{{\bf y}}

\def\f{\mbox{field}}
\def\a{\alpha}
\def\lx{\xi_x}
\def\ly{\eta_y}

\def\PML{{\rm \text{PML}}}
\def\B{B}
\def\P{\mathcal{P}}

\def\L{{\mathcal{L}}}

\def\Nbx{N_1}
\def\Nby{N_2}

\def\half{\frac{1}{2}}
\def\f32{\frac{3}{2}}
\def\epml{e^{-\frac{1}{2} k \gamma_{0} \bar{\sigma} }}

\def\s1{{s-1}}
\def\s2{{s-2}}

\def\nSolveDDM{n_{\text{DDM\textvisiblespace Solv}}}
\def\nIterDDM{n_{\text{DDM\textvisiblespace Iter}}}
\def\nIterGMERES{n_{\text{GMRES\textvisiblespace Iter}}}
\def\nLocalSolve{n_{\text{Local\textvisiblespace Solv}}}

\begin{document}

\newtheorem{theorem}{Theorem}[section]
\newtheorem{lemma}[theorem]{Lemma} 
\newtheorem{remark}[theorem]{Remark}

%
\newtheorem{algorithm_}{Algorithm}[section]

\title{An Additive Overlapping Domain Decomposition Method 
 for the Helmholtz Equation}

\author{Wei Leng \thanks{State Key Laboratory of Scientific and Engineering Computing,
Chinese Academy of Sciences, Beijing 100190, China. Email: {\tt wleng@lsec.cc.ac.cn}. W. Leng's research is partially supported by Natural Science Foundation of China under grant number 11501553, and the National Center for
Mathematics and Interdisciplinary Sciences of  Chinese Academy of
Sciences.}\and
Lili Ju \thanks{Department of Mathematics, University of South Carolina,
Columbia, SC 29208, USA. Email: {\tt ju@math.sc.edu}. L. Ju's research is partially supported by US National Science Foundation under grant number
DMS-1521965.}
}

\maketitle

\begin{abstract}
In this paper, we propose and analyze an additive domain decomposition method (DDM) 
for solving the high-frequency Helmholtz equation with the Sommerfeld radiation condition. In the proposed method, the computational domain is partitioned into structured subdomains along all spatial directions, and each subdomain contains an overlapping region for source transferring.   At each iteration all subdomain PML problems are solved  completely in parallel, then all horizontal, vertical and corner directional residuals on each subdomain are passed to its corresponding neighbor subdomains 
as the source for the next iteration. This  DDM method is highly scalable in nature and theoretically shown to produce the  exact solution for the PML problem defined in $\R^2$ in the constant medium case.  A slightly modified version of the method for bounded truncated domains is also developed for its use in practice  and an  error estimate is  rigorously proved.
Various numerical experiments  in two and three dimensions are conducted on the  supercomputer ``Tianhe-2 Cluster'' to verify the theoretical results and demonstrate excellent performance of  the proposed method as an iterative solver or a preconditioner. 
\end{abstract}

\textbf{keywords}

Additive domain decomposition method, Helmholtz equation, perfectly matching layer, convergence analysis, parallel computing


\pagestyle{myheadings}
\thispagestyle{plain}
\markboth{Wei Leng and Lili Ju}{Additive overlapping DDM for Helmholtz Equation}


\section{Introduction}

In this paper,  we consider the well-known Helmholtz equation in the space $\R^n$
\begin{align} \label{eq:helm}
  \Delta u + k^2 u &= f ,\qquad  \mbox{in} \,\,\, \R^n
\end{align}
with the Sommerfeld radiation condition
\begin{align} \label{cond_s}
  r^{\frac{n-1}{2}} (\frac{\p u}{\p r} - \mathbf{i} k u) &\rightarrow 0, \qquad \mbox{as} \,\,\, r = |\bx| \rightarrow \infty, 
\end{align}
where $u(\bx)$ is the unknown function, $f(\bx)$ is the source, $k(\bx)$ denotes the wave number, and $\DD k(\bx) = {\omega}/{c(\bx)}$,
where $\omega$ is the angular frequency and $c(\bx)$ is the wave speed.
Among many discretization and solution techniques for solving the Helmholtz equation \eqref{eq:helm}, 
the domain decomposition method (DDM) \cite{Schwarz1870, Lions1989} is a very effective  approach when the problem scale is very large  and has been thoroughly studied for 
several decades. 
The basic idea is as follows: truncated with the absorbing boundary condition,  the subdomain problem is (exactly or approximately) solved during each iteration, and  the subdomain solutions is  then passed to their respective neighbor 
subdomains via  interface conditions to carry on the wave propagation process.   The DDM is also especially popular to be  used as 
a preconditioner to the iterative solver for the global discrete system.

The DDM was first used to solve the Helmholtz problem by Despr\'{e}s
in \cite{Despres1990}, and then various DDM algorithm has been
developed to solve the Helmholtz problem
\cite{Schwarz1,Schwarz2,Schwarz3,Schwarz4,Schwarz5,Schwarz6}, these 
algorithm are categorized as the Schwarz methods with or without overlap. The boundary conditions
at the subdomain interfaces affects the convergence of rate of the
Schwarz method, thus Gander introduced the optimized Schwarz methods
using an optimal non-local boundary condition
\cite{SchwarzBdry1,SchwarzBdry2,SchwarzBdry3,SchwarzBdry4,SchwarzBdry5,SchwarzBdry6}.
A non-overlapping Schwarz method was proposed by Boubendir et al. \cite{Boubendir2012} using Pad\'{e} approximations of the Dirichlet to Neumann (DtN) map,
and the effective convergence is quasi-optimal.

The DDM with the transmission condition that involves sources yield fast methods to solve
the high-frequency Helmholtz equation, the sweeping DDM preconditioner proposed by
Engquist and Ying \cite{Engquist2011a, Engquist2011b} was the first of this type, and is shown to be very effective,
and followed by the development of 
various sweeping type DDM, which differs mainly at the transmission condition at the interfaces of subdomains.
The transmission condition of Engquist and Ying \cite{Engquist2011a, Engquist2011b} is that the residual in just one layer is taken as the source.
The DDM  introduced by Stolk \cite{Stolk2013} uses a transmission condition that 
the derivative of the solution is treated as a delta source.  
The double sweeping preconditioner developed by Vion an Geuzaine \cite{Vion2014} uses a mixed boundary condition
that involves DtoN map as the transmission condition.
The polarized trace method developed by Zepeda \cite{Zepeda2014} impose a transmission condition
that both single and double layer potentials are taken as delta sources.  
Instead of using the delta source approach, the source transfer DDM (STDDM) \cite{Chen2013a,Chen2013b} proposed by Chen and Xiang
uses a smooth source coming from the residual as the transmission condition, 
and furthermore, an error estimate of the method has been rigorously proven. 
It is noted that the sweeping type DDMs often partition the domain in one direction, and the sweeping 
order would be along the direction forwards and backwards.
Two orthogonal directions sweeping could be done in a recursive way, for instance,  as proposed by Liu \cite{Liu2015b}
using the  sweeping preconditioner \cite{Engquist2011a, Engquist2011b},
and by Wu \cite{Wu2015} using the source transfer DDM \cite{Chen2013a,Chen2013b}.
The sweeping type DDM is {\em multiplicative}, 
the subdomain problems are solved one after another in a particular order/direction, 
which could be interpreted as the process of $LU$ or $LDL^T$ factorization. 
The multiplicative DDMs usually have better convergence  than the additive DDMs,
however, the sequential sweeping order causes many  difficulties to parallel computing in term of efficiency and scalability.
In addition, the subdomain problems in the DDMs are often solved with  direct solvers, such as 
the multi-frontal method \cite{Frontal1, Frontal2, MUMPS},
or the the multi-frontal method using hierarchically semiseparable structure
\cite{HSS, Wang2016},
but such one-direction partitions of the domain  leave us  large subdomains in practice and 
 the  robustness of the direct solvers is still a challenge for such large size problems.
 

 There are not many research on the of additive DDM
  with the transmission condition that involves sources for solving high-frequency Helmholtz equation. 
 Liu and Ying proposed an additive sweeping
 preconditioner in \cite{Liu2015a}, where the computational domain is
 partitioned in one direction into many thin layers, and the
 subproblems on all respective layers are solved in parallel with
 direct solvers.  Wei introduced an additive overlapping DDM
 solver in \cite{Leng2015}, which uses a modified source transfer
 technique. The main advantages of the additive DDMs is that they are
 more suitable and scalable for parallel computing and much easier to
 implement on massively distributed machines than multiplicative ones.
In this paper, we propose, analyze and test a new additive overlapping DDM for the Helmholtz equation \eqref{eq:helm}
based on the source transfer technique \cite{Chen2013a}. We focus on illustrations of the DDM method 
on two-dimensional problems with constant medium, but it is very natural to use the method as an iterative solver or a preconditioner  for the variable medium problems, and it is very straightforward to extend the proposed algorithm to three-dimensional problems.  
In our method,   all the subdomain PML problems are solved completely in parallel at each iteration,
and all {\em horizontal}, {\em vertical} and {\em corner} directional residuals on each subdomain
are then passed to its corresponding neighbor subdomains as the source for the next iteration.
Based on the source transfer analysis, the DDM method for the PML problem defined in $\R^2$ is theoretically shown to produce the  exact solution in the  constant medium case.
Since the practical computations of the PML problems  only can be done in bounded truncated domains,
we then slightly modify the  additive overlapping DDM to be useful for the truncated PML problem (which is  an approximation to the full space problem).
An error estimate  of the method is further derived, and 
in particular, our proof revises and simplifies
some analysis techniques used in  \cite{Chen2013a} and successfully generalizes them to the additive type of DDM with the concurrent residual transfers in all directions.

The rest of the paper is organized as follows. In Section 2, we first present the  
additive overlapping  DDM for the PML problem in the full space of $\R^2$, and show that the DDM solution is
exactly the solution to the problem  in the  constant medium case.  In Section 3,
the corresponding   DDM in bounded truncated domains is then developed  for  practical use 
and its error estimate is also  derived.
In Section 4, various numerical experiments in two-dimensional  and three-dimensional spaces are performed to verify the
 theoretical results and demonstrate effectiveness and scalability of the DDM method as an iterative solver or a preconditioner
for the  global discrete system. 
Some concluding remarks are finally drawn in Section 5.

\section{The additive overlapping DDM in $\R^2$}

In this section, the PML method and the source transfer technique are first briefly reviewed, then we will develop and analyze an additive overlapping DDM for solving  the PML equation in the whole space  $\R^2$. The proposed DDM method  will make use of  
the structured overlapping domain decomposition together with the source transfer in  horizontal, vertical and  corner directions.  
The constant medium problem (the  constant wave number $k(\bx)\equiv k$) is assumed for development and analysis of our DDM method.

\subsection{Perfect match layer and source transfer}

The solution of the Helmholtz equation \eqref{eq:helm}  with the Sommerfeld radiation condition  \eqref{cond_s}  in the constant medium case
can be represented by
\begin{equation}
u(\bx) = - \int_{\R^2} f(\by) G(\bx, \by) \;d\by,\qquad \forall\,\bx\in\R^2,
\end{equation}
where $\DD G(\bx,\by) = \frac{\mathbf{i}}{4} H_0^{(1)} (k|\bx-\by|)$ is the fundamental solution of the
Helmholtz equation defined by
\begin{equation}
\Delta G(\bx,\by) + k^2 G(\bx,\by) = -\delta_{\by}(\bx),\qquad\forall\,\bx\in\R^2,
\end{equation}
Note that $H_0^{(1)}(z)$, $z \in \mathbb{C}$ is often referred as the first Hankel function of order zero \cite{Ar85}.

Consider a  rectangular box
 $\B = \{(x_1,x_2)^T: a_j \le x_j \le b_j, j = 1,2\}\in\R^2$, 
and the center of the rectangle is denoted by $(c_1, c_2)$, where $c_j = \frac{a_j+b_j}{2}$, $j = 1,2$.
We will build and solve the PML equation in $\R^2$ for the rectangular box $\B$, for instance,
 replace the solution outside  $\B$  by 
the complex coordinate stretching as done in the so-called uniaxial PML method \cite{Chew1994, Kim2010, Bramble2013, Chen2013a}.
Let $\alpha_1(x_1) = 1 + \ii \sigma_1(x_1) $,  $\alpha_2(x_2) = 1 + \ii \sigma_2(x_2) $,
with $\sigma_1,\sigma_2$  being  piecewise smooth functions such that
\begin{equation}
  \sigma_j(x) = \left\{  
      \begin{array}{ll}
        \tilde{\sigma}(x_j - b_j), & \text{if} \,\,\, b_j \le x_j,\\
        0,                         & \text{if} \,\,\, a_j < x_j < b_j,\\
        \tilde{\sigma}(a_j - x_j), & \text{if} \,\,\, x_j \le a_j,\\
      \end{array}
  \right. 
\end{equation}
where $\tilde{\sigma}(t)$ is a smooth medium profile function satisfying that $\tilde{\sigma}(t) = \gamma_0$, for $t > d$, 
and both $d$ and $\gamma_0$ are some positive constants.
Then the complex coordinate stretching $\tilde{\bx}(\bx)$ for $\bx = (x_1, x_2)$ is defined as
\begin{equation}
\tilde{x}_j(x_j) = c_j+\int_{c_j}^{x_j} \alpha_j(t) dt = x_j
          + \ii \int_{c_j}^{x_j} \sigma_j(t) dt,\qquad j = 1,2.
\end{equation}
Denote $z^{1/2}$ the analytic branch of $\sqrt{z}$ such that $\mbox{Re}(z^{1/2}) > 0$
 for $z \in \mathbb{C}[0,+\infty]$, the distance in the complex plane is defined as
\begin{equation}
\rho(\tilde{\bx}, \tilde{\by}) = \big[(\tilde{x}_1(x_1) - \tilde{y}_1(y_1))^2 +  (\tilde{x}_2(x_2) - \tilde{y}_2(y_2))^2 \big]^{1/2}.
\end{equation}  
Now we define 
\begin{equation}
\tilde{u}(\bx) = u(\tilde{\bx}(\bx)) = -\int_{\R^2} f(\by) G(\tilde{\bx}, \tilde{\by}) d\by, \qquad \forall\, \bx \in \R^2,
\end{equation}
where $\DD G(\tilde{\bx}, \tilde{\by}) = \frac{\mathbf{i}}{4} H_0^{(1)} (k \rho(\tilde{\bx}, \tilde{\by}) )$. 
We also assume that  $f$ is compactly supported in $\B$,  
then $\tilde{u}$ is well-defined in $H^1_{loc}(\R^2)$,  satisfies $\tilde{u}(\bx) = u(\bx)$ in $\B$, and
decays exponentially as $|x|\rightarrow \infty$. Consequently,
  it is the solution of the following  PML equation with the source $f$:
\begin{equation} \label{eq:PML}
J_{\B}^{-1} \nabla \cdot (A_{\B} \nabla\tilde{u}) + k^2 \tilde{u} = f, \qquad\forall\,\bx\in\R^2,
\end{equation}
where $\DD A_{\B}(x) = \mbox{diag}\left(\frac{\a_2(x_2)}{\a_1(x_1)}, \frac{\a_1(x_1)}{\a_2(x_2)}\right)$ and $J_{\B}(x) = \a_1(x_1) \a_2(x_2)$.
We denote by $\L_{\B}:=J_{\B}^{-1} \nabla \cdot (A_{\B} \nabla  \, \boldsymbol{\cdot}) + k^2  $  the linear operator associated with \eqref{eq:PML}.
The weak formulation of the PML problem \eqref{eq:PML} can be given as follows: for $f \in H^1(\R^2)^{\prime}$,
find $u \in H^1(\R^2)$ such that
\begin{equation}
  \label{eq:PMLweak}
  (A_{\B} \nabla \tilde{u}, v) - k^2 (J_{\B}\tilde{u}, v) = - \langle J_{\B}f, v \rangle, \qquad \forall\, v\in H^1(\R^2).
\end{equation}
where $(\cdot,\cdot)$ denotes the inner product in $L(\R^2)$ and $\langle \cdot \rangle$  the duality pairing 
between $H^1(\R^2)^{\prime}$ and $H^1(\R^2)$.
The well-posedness of the problem \eqref{eq:PMLweak} has been established and could be found in  \cite[Lemma 3.3]{Chen2013a}. Let us denote by $\P_{\B}$ the above PML problem \eqref{eq:PMLweak}
associated with the rectangular box $\B$.

\begin{figure}[ht!]	
	\begin{center}
		\includegraphics[width=0.6\textwidth]{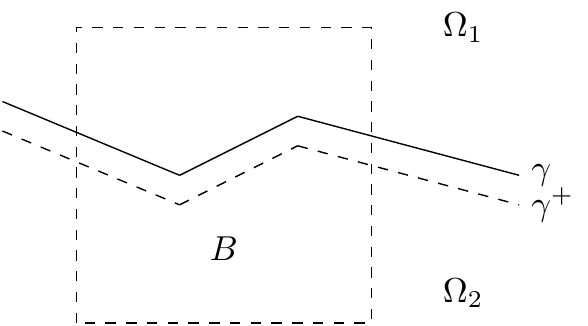}
		\caption{Divide $\R^2$ and the box $B$ with  $\gamma$.}\label{fig:twopart}
	\end{center}
\end{figure}

The source transfer technique in $\R^2$  is briefly restated as follows. 
Suppose that $\R^2$ is divided into two parts $\Omega_1$ and $\Omega_2$ 
by the piecewise smooth curve $\gamma$, and at the same time  $\gamma$ also divides the rectangular box $\B$ into two parts.
Let $\Omega_1^+ = \{\bx : \rho(\bx, \Omega_1) \le \tilde{d}\}$, where $\tilde{d}$ is a constant,
 and $\gamma^+ = \p \Omega_1^+$,
as show in Figure \ref{fig:twopart}. There always exists
a smooth function $\beta \in C^2(\R^2)$  such that 
$$\beta|_{\Omega_1} \equiv 1, \quad \beta|_{\R^2 \setminus \Omega_1^+} \equiv 0, \quad
0 \leq \beta \leq 1,$$ 
and 
$$ |\nabla \beta(\bx) | < C, \qquad \forall\,\bx\in \Omega_1^+\setminus\Omega_1,$$
where $C$ denotes some generic positive constant.

It is quite straightforward to prove the following result. Denote by ${\chi}_{\Omega_2}$ the characteristic function of $\Omega_2$.
\begin{lemma} \label{thm:main}
Suppose the support of $f$ is in $\Omega_1 \cap \B $.
Let $u_0$ be the solution to the PML problem $\P_B$ with  the source $f$ (i.e, $\L_{\B} u_0 = f$ in $\R^2$).
Given $u_1 \in H^1(\R^2)$ such that $u_1 = u_0 \beta$ in $\R^2$ and let $u_2$ be solution to 
 the PML problem $\P_B$ with  the source $-(\L_{\B} u_1)\chi_{\Omega_2}$  (i.e.,
$	\L_{\B} u_2 = -(\L_{\B} u_1)\chi_{\Omega_2}$ in $\R^2$). Then we have that $u_1 + u_2 = u_0$ in $\R^2$ and  $u_2 = 0$ in $\Omega_1$.
\end{lemma}

The two domains in Lemma \ref{thm:main} have an overlapping region
$\Omega_1^+ \setminus \Omega_1$.   
For convenience,
we shift the media profile function $\tilde{\sigma}(t)$ by $\tilde{d}$ so that 
\begin{equation} \label{eq:PMLprofile}
  \tilde{\tilde{\sigma}}(t) = \left\{  
      \begin{array}{ll}
        0,                          & \text{if} \,\,\, t \leq \tilde{d}\\
        \tilde{\sigma}(t-\tilde{d}), & \text{if} \,\,\, t > \tilde{d}\\
      \end{array}
  \right. 
\end{equation}
where $\tilde{\tilde{\sigma}}(t)$ is the shifted medium profile.
From now on, we will denote $\tilde{\sigma}(t)$ as the above shifted medium profile for simplicity, 
and in this way, for the PML problem $\P_{\B}$, an extended region $\{\bx: \rho(\bx, \B) \leq \tilde{d} \} \setminus \B$ is always reserved as a possible overlapping region.

\begin{figure}[ht!]	
	\begin{center}
		\includegraphics[width=.4\textwidth]{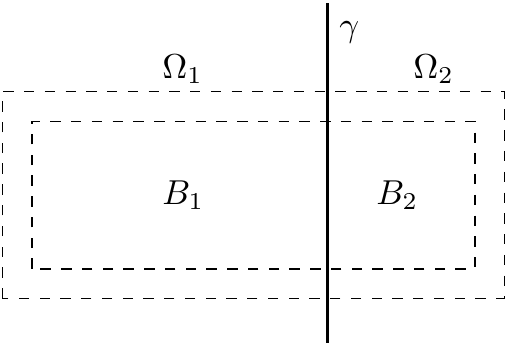}
		\hspace{1cm}
		\includegraphics[width=.4\textwidth]{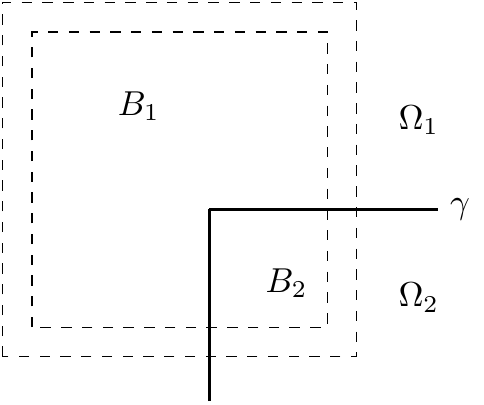}
		\caption{Source transfer in horizontal direction (left) and in the lower-right corner direction (right).}\label{fig:transfer}
	\end{center}
\end{figure}

We  will propose two types of source transfers in our method  for the two-dimensional problem.
Let $\kappa_1$ and $\kappa_2$ be two arbitrary constants such that $a_j \leq \kappa_j \leq b_j$, $j=1,2$.
The first type is the transfer in $x$- or $y$-direction,
as is shown in Figure \ref{fig:transfer} (left), where $\gamma = \{(x_1, x_2): x_1 = \kappa_1\}$ 
or $\gamma = \{(x_1, x_2): x_2 = \kappa_2\}$. We remark that such type of transfer is  similar to the one used in STDDM \cite{Chen2013a}.
The second type is the transfer in any of four corner directions,
as is shown in Figure \ref{fig:transfer} (right), where 
$\gamma = \{(x_1, x_2): x_1 = \kappa_1, x_2 \geq \kappa_2\} \cup \{(x_1, x_2): x_1 \leq \kappa_1, x_2 = \kappa_2\}$.
In both cases, we have $u_2 = 0$ in $\Omega_1$, therefore instead of solving the PML problem $\P_B$,
we can turn to solve the PML problem $\P_{B \cap \Omega_2}$ whose solution is zero in $\Omega_1$. 

\subsection{The  DDM for the PML problem in $\R^2$}

We now develop an additive overlapping DDM method to solve
the PML equation \eqref{eq:PMLweak}  associated with a rectangular domain $\Omega = (-l_1,l_1)\times(-l_2,l_2)$ in $\R^2$.
 Assume that  $\Omega$ is uniformly partitioned into
 $\Nbx \times \Nby$ non-overlapping rectangular subdomains. 
 Let $\Delta\xi = 2 l_1 / \Nbx$,  
 $\xi_i = -l_1 + (i-1) \Delta \xi$, $i = 1, \ldots, \Nbx+1$, 
 and $\Delta\eta = 2 \l_2 / \Nby$, 
 $\eta_j = -l_2 + (j-1) \Delta\eta$, $j = 1, \ldots, \Nby+1$.
 Then we have $\Nbx\times \Nby$ non-overlapping rectangular subdomains as
 $$\Omega_{i,j}
 : = [\xi_{i}, \xi_{i+1}] \times [\eta_{j}, \eta_{j+1}],\qquad i = 1, 2,\ldots, \Nbx, j = 1, 2,\ldots, \Nby.$$
It is clear that the PML equation associated with each  rectangular subdomain $\Omega_{i,j}$ is needed to be solved in the DDM method.
The source $f$, which is assumed to be compactly supported in  $\Omega$, is decomposed to
$$f_{i,j} = f \cdot \chi_{\Omega_{i,j}} , \qquad i = 1, \ldots, \Nbx, j = 1, \ldots, \Nby.$$
For convenience, we also define the box $\Omega_{i0,i1;j0,j1}$ ($1\leq i_0<i_1\leq N_1+1$, $1\leq j_0<j_1\leq N_2+1$), which consists of a set of  rectangular subdomains as follows
$$\Omega_{i0,i1;j0,j1} := \bigcup\limits_{i0 \leq i \leq i1 \atop j0 \leq j \leq j1} \Omega_{i, j}.$$
Notice that the PML profile as \eqref{eq:PMLprofile} makes
each subdomain has an overlapping region with its neighbor subdomains, thus we next define an overlapping domain decomposition of the two-dimensional space $\R^2$ as 
 $$\widetilde{\Omega}_{i,j}
 : = (\tilde{\xi}_{i}-\tilde{d}, \tilde{\xi}_{i+1}+\tilde{d}) 
   \times (\tilde{\eta}_{j} - \tilde{d}, \tilde{\eta}_{j+1} + \tilde{d}),
   \quad i = 1, \ldots, \Nbx, j = 1, \ldots, \Nby,$$
where 
\begin{align*}
\tilde{\xi}_{i} &= \left\{ 
\begin{array}{ll}
-\infty,  &\,\,\, i = 1,\\
\xi_i,  & \,\,\, i = 2,\ldots,\Nbx,\\
+\infty, & \,\,\, i = \Nbx+1,
\end{array}
\right. &  
\tilde{\eta}_{j} &= \left\{
\begin{array}{ll}
-\infty,  &\,\,\, j = 1,\\
\eta_j,  & \,\,\, j = 2,\ldots,\Nby,\\
+\infty, & \,\,\, j = \Nby+1.
\end{array}
\right. 
\end{align*}

A few notations are introduced as follows.  
Let
$\beta_0(t)$ be a monotone function in $C^2(\R)$, such that 
$\beta_0(t) = 1$ for $t \leq 0$, 
$\beta_0(t) = 0$ for $t \geq 1$,
and $|\beta_0^{\prime}(t)| < C$ for $ 0 < t < 1.$ 
Then define
\begin{align*}
\beta_{\leftarrow; i} &= \left\{ 
      \begin{array}{ll}
      1,  &\, i = 1,\\
      \beta_0(\frac{\xi_i - x_1}{\tilde{d}}),  & \, i = 2,\ldots,\Nbx+1,
      \end{array}
\right. &  \beta_{\downarrow; j} &= \left\{
      \begin{array}{ll}
      1,  &\, j = 1,\\
      \beta_0(\frac{\eta_j - x_2}{\tilde{d}}),  & \, j = 2,\ldots,\Nby+1,
      \end{array}
\right. \\
\beta_{\rightarrow; i} &= \left\{
      \begin{array}{ll}
      \beta_0(\frac{x_1 - \xi_i}{\tilde{d}})  & \, i = 1,\ldots,\Nbx, \\
      1,  &\, i = \Nbx+1,
      \end{array} 
\right. &  \beta_{\uparrow; j} &= \left\{ 
      \begin{array}{ll}
      \beta_0(\frac{x_2 - \eta_j}{\tilde{d}})  & \, j = 1,\ldots,\Nby, \\
      1,  &\, j = \Nby+1,
      \end{array}
\right. 
\end{align*}

 and 
 \begin{align*}
 \beta_{\swarrow; i,j} &= \beta_{\leftarrow; i} \beta_{\downarrow;j} \,\, ,&
 \beta_{\searrow; i,j} &= \beta_{\rightarrow; i} \beta_{\downarrow;j}  \,\, , &\\
 \beta_{\nwarrow; i,j} &= \beta_{\leftarrow; i} \beta_{\uparrow;j}  \,\,,&
 \beta_{\nearrow; i,j} &= \beta_{\rightarrow; i} \beta_{\uparrow;j}  \,\,, 
 & \beta_{0; i,j} &=  \beta_{\leftarrow; i} \beta_{\rightarrow;i+1} \beta_{\downarrow; j} \beta_{\uparrow;j+1} .
 \end{align*}
These functions  will be used as $\beta$ in Lemma \ref{thm:main} 
for  different cases.
Denote $\L_{i,j}$ as the linear operator associated with the PML problem $\P_{\Omega_{i,j}}$.
Define the following characteristic functions for the half spaces and the quarter spaces: 
\begin{align*}
 \chi_{\leftarrow; i}   &:= \chi_{(-\infty, \xi_i)\times(-\infty,+\infty)}, 
&\chi_{\rightarrow; i}  &:= \chi_{(\xi_{i}, +\infty) \times (-\infty,+\infty)}, \\
 \chi_{\downarrow; j}   &:= \chi_{(-\infty,+\infty) \times (-\infty, \eta_j)}, 
&\chi_{\uparrow; j}     &:= \chi_{(-\infty,+\infty) \times (\eta_{j}, +\infty) }, \\
  \chi_{\swarrow; i,j}  &:=  \chi_{(-\infty, \xi_i)\times(-\infty, \eta_j)}, 
 &\chi_{\searrow; i,j}  &:=   \chi_{(\xi_{i}, +\infty)\times(-\infty, \eta_j)}, \\
  \chi_{\nwarrow; i,j}  &:=  \chi_{(-\infty, \xi_i)\times(\eta_{j}, +\infty)}, 
 &\chi_{\nearrow; i,j}  &:=  \chi_{(\xi_{i}, +\infty)\times(\eta_{j}, +\infty)}.
\end{align*}

Using the $\beta$ functions, characteristic functions and linear operators  defined above,
we are able to define the transfer operator as follows:
\begin{align*}
 \Psi_{\leftarrow; i,j} (v) &:= -\L_{i-1,j} (\beta_{\leftarrow;i} v)  \chi_{\leftarrow; i-1},\\
\Psi_{\rightarrow; i,j}( v) &:= -\L_{i+1,j} (\beta_{\rightarrow;i+1} v) \chi_{\rightarrow; i+1}, \\
 \Psi_{\downarrow; i,j} (v) &:= -\L_{i,j-1} (\beta_{\downarrow;j} v)  \chi_{\downarrow; j-1},\\
\Psi_{\uparrow; i,j} (v)   &:= -\L_{i,j+1} (\beta_{\uparrow;j+1} v) \chi_{\uparrow; j+1},\\
 \Psi_{\swarrow; i,j} (v)   &:= -\L_{i-1,j-1} (\beta_{\swarrow;i,j} v) \chi_{\swarrow; i-1,j-1}, \\
\Psi_{\searrow; i,j} (v)   &:= -\L_{i+1,j-1} (\beta_{\searrow;i+1,j} v) \chi_{\searrow; i+1,j-1},\\ 
\Psi_{\nwarrow; i,j} (v)   &:= -\L_{i-1,j+1} (\beta_{\nwarrow;i,j+1} v) \chi_{\nwarrow; i-1,j+1},\\
\Psi_{\nearrow; i,j} (v )  &:= -\L_{i+1,j+1} (\beta_{\nearrow;i+1,j+1} v) \chi_{\nearrow; i+1,j+1}.
\end{align*}
Note that there are a total of $2^3-1=8$ directions for source transfer since each subdomain could have up to 8 neighbor subdomains
in the two-dimensional space.
%
%

%
%
%
%

\begin{figure}[ht!]	
	\begin{center}
		\begin{minipage}{.4\textwidth}
		\includegraphics[width=\textwidth]{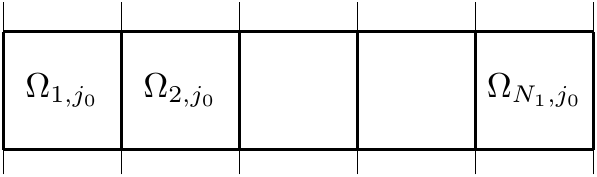}
		\end{minipage}
		\hspace{0.5cm}
	    \begin{minipage}{.45\textwidth}
		\includegraphics[width=\textwidth]{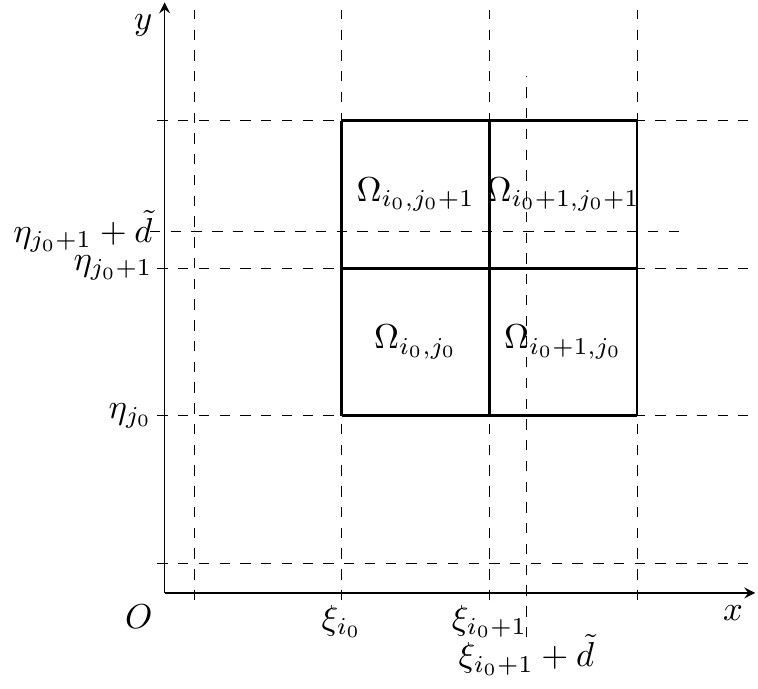}
		\end{minipage}
		\caption{Source transfers in the horizontal direction (left) and in the upper-left corner direction (right) in the DDM.} \label{fig:transf_dir1}
	\end{center}
\end{figure}

First let us illustrate the source transfer in the {\em horizontal} or {\em vertical} direction of the domain decomposition.  
The procedure is essentially the same as the STDDM \cite{Chen2013a}.
Assume $f_{1,j_0} \neq 0$, the PML problem 
$\P_{\Omega_{1,\Nbx;j_0,j_0}}$ 
with the source $f_{1,j_0}$ is to be solved using the DDM, and its exact solution is denoted by $u$
(see Figure \ref{fig:transf_dir1} (left) for illustration).
The subdomain PML problem on $\Omega_{i, j_0}$ can be  solved in order from $i = 1$ to $i = \Nbx$
to construct the  solution $u$:
first, the problem $\P_{\Omega_{1, j_0}}$ with the source $f_{1,j_0}$ is solved, the solution is denoted as $u_1$,
we have $u_1 = u$ in $\{(x_1,x_2): x_1 \leq \xi_2 + \tilde{d}\}$; 
second, solve the PML problem $\P_{\Omega_{2, j_0}}$ with
the source $\Psi_{\rightarrow;1,j_0} (u_1) = -\L_{{2,j_0}} (\beta_{\rightarrow; 2} u_1) \chi_{\rightarrow; 2}$, the solution is denoted $u_2$,
then using Lemma \ref{thm:main} on the box $\Omega_{1,2;, j_0, j_0}$, we know that
$u_1 \beta_{\rightarrow;2} + u_2 = u$ in $\{(x_1,x_2): x_1 \leq \xi_3 + \tilde{d}\}$; 
third, solve the PML problem $\P_{\Omega_{3, j_0}}$ with
the source $\Psi_{\rightarrow;2,j_0} (u_2) = -\L_{{3,j_0}} (\beta_{\rightarrow; 3} u_2) \chi_{\rightarrow; 3} $, the solution is denoted $u_3$,
then using Lemma \ref{thm:main} on box $B_{1,3;, j_0, j_0}$, we know that
$u_1 \beta_{\rightarrow;2} + u_2\beta_{\rightarrow;3} + u_3= u$ in $\{(x_1,x_2): x_1 \leq \xi_4 + \tilde{d}\}$; 
repeat the process until the last Box $\Omega_{\Nbx, j_0}$,
and the DDM solution is found as $\sum\limits_{i=1,\ldots,\Nbx} u_i\beta_{\rightarrow;i+1}= u.$

\begin{remark}
When the source transfer is in the horizontal or vertical direction, the source of the subdomain
problem possesses a specific form. Let 
\begin{align} \label{eq:Phi}
\Phi_{m; i,j} (u; \beta) := J_{\Omega_{i,j}}^{-1} \frac{\p }{\p x_m} 
     \big( (A_{\Omega_{i,j}})_{m,m} \frac{\p \beta}{\p x_m} u  \big)  
     + J_{\Omega_{i,j}}^{-1} \frac{\p \beta}{\p x_m} 
      (A_{\Omega_{i,j}})_{m,m} \frac{\p u}{\p x_m},
\end{align}
where $m = 1,2$,
then the above  sources  for the subdomain problem could be represented respectively  as   
\begin{align} \label{eq:transf_form_x}
\Psi_{\rightarrow ; i,j_0}(u_i) = \Phi_{1; i+1,j_0} (u_i; \beta_{\rightarrow; i+1}) \chi_{\rightarrow;i+1}, \qquad i = 1, 2, \dots, \Nbx-1,
\end{align}
where the solution of other subdomain problem and its first order derivatives are involved but the wave number $k$ is not. 
Similar forms exist for other horizontal or vertical directional source transfers.  
\end{remark}

%
%
%
%
Next we  illustrate the source transfer in the corner directions of the domain decomposition.  
Suppose $f_{i_0,j_0} \neq 0$, 
the PML problem $\P_{\Omega_{i_0,i_0+1;j_0,j_0+1}}$ with the source $f_{i_0,j_0}$ is to be solved with the DDM,
and its solution is denoted as $u$ again (see Figure \ref{fig:transf_dir1} (right) for illustration).
First the solution to the PML problem $\P_{\Omega_{i_0,j_0}}$ with the source $f_{i_0,j_0}$ is solved, 
and the solution is denoted as $u_0$,
 obviously,
\begin{equation} \label{eq:corner_u1}
u_0 = u, \,\,\, \text{in} \,\,\, \{(x_1,x_2): x_1 \leq \xi_{i_0+1} + \tilde{d} \,\,\,\, \text{and}\,\,\,\, x_2 \leq \eta_{j_0+1} + \tilde{d} \}.
\end{equation} 
Then the horizontal source transferring is applied, 
the problem $\P_{\Omega_{i_0+1,j_0}}$ with the source $\Psi_{\rightarrow; i_0,j_0}(u_0)$ is solved, 
and the solution is denoted as $u_{\rightarrow}$, we have
\begin{align}
u_{0} \beta_{\rightarrow;\xi_{i_0+1}} + u_{\rightarrow} &= u, \qquad  \text{in}  \,\,\, \{(x_1,x_2): x_2 \leq \eta_{j_0} + \tilde{d}\}, \label{eq:corner_u2}\\
u_{\rightarrow} &= 0, \qquad \text{in} \,\,\, \{(x_1,x_2): x_1 \leq \xi_{i_0} \}.  \label{eq:corner_u22}
\end{align} 
Similarly, the vertical source transferring is applied, 
the PML problem $\P_{\Omega_{i_0,j_0+1}}$ with the source $\Psi_{\uparrow; i_0,j_0}(u_0)$ is solved,
and the solution is denoted as $u_{\uparrow}$, we have 
\begin{align}
  u_0 \beta_{\uparrow;\eta_{j_0+1}} + u_{\uparrow} &= u, \qquad \text{in} \,\,\, \{(x_1,x_2): x_1 \leq \xi_{i_0} + \tilde{d}\}, \label{eq:corner_u3} \\
  u_{\uparrow} &= 0, \qquad \text{in}\,\,\,  \{(x_1,x_2): x_2 \leq \eta_{j_0} \}. \label{eq:corner_u33}
\end{align}
By \eqref{eq:corner_u1}-\eqref{eq:corner_u33}, we obtain that 
\begin{equation} \label{eq:corner_u}
u_0 \beta_{\nearrow;\xi_{i_0+1,j_0+1}} 
+ u_{\rightarrow} \beta_{\uparrow;\eta_{j_0+1}}
+ u_{\uparrow}  \beta_{\rightarrow;\xi_{i_0+1}}
= u (\beta_{\rightarrow;\xi_{i_0+1}}
	+ \beta_{\uparrow;\eta_{j_0+1}}
	- \beta_{\nearrow;\xi_{i_0+1,j_0+1}} ),  \;\;\text{in} \,\, \R^2.
\end{equation}
At last the PML problem $\P_{\Omega_{i_0+1, j_0+1}}$  with the source
\begin{align}
 &-\L_{B_{i_0,i_0+1;j_0,j_0+1}} \big(u (\beta_{\rightarrow;\xi_{i_0+1}}
 + \beta_{\uparrow;\eta_{j_0+1}}
 - \beta_{\nearrow;\xi_{i_0+1,j_0+1}} ) 
 \big)  \chi_{\nearrow; i_0+1,j_0+1}, \nonumber \\
&\quad=\;-\L_{i_0+1;j_0+1} \big(u_0 \beta_{\nearrow;\xi_{i_0+1,j_0+1}} 
       + u_{\rightarrow} \beta_{\uparrow;\eta_{j_0+1}}
       + u_{\uparrow} \beta_{\rightarrow;\xi_{i_0+1}} \big) \chi_{\nearrow; i_0+1,j_0+1} \nonumber\\
&\quad=\; \Big( \Psi_{\nearrow; i_0,j_0}(u_0) 
       + \Psi_{\uparrow; i_0+1,j_0}(u_{\rightarrow})
       + \Psi_{\rightarrow; i_0,j_0+1}(u_{\uparrow})\Big) \chi_{\nearrow; i_0+1,j_0+1}, \label{eq:rhs}
\end{align}
is solved and the solution is denoted as $u_{\nearrow}$.
Observe that
$$
\begin{array}{l}
\beta_{\rightarrow;\xi_{i_0+1}}
+ \beta_{\uparrow;\eta_{j_0+1}}
- \beta_{\nearrow;\xi_{i_0+1,j_0+1}} 
= 1 - (1-\beta_{\rightarrow;\xi_{i_0+1}})(1-\beta_{\uparrow;\eta_{j_0+1}})\\
\qquad=\left\{
\begin{array}{ll}
1 & {\rm in} \;\;\{(x_1,x_2):  x_1 \leq \xi_{i_0+1} \,\,\,  \text{or} \,\,\, x_2 \leq \eta_{j_0+1}\},\\
0 & {\rm in} \;\; \{(x_1,x_2):  x_1 \geq \xi_{i_0+1} + \tilde{d} \,\,\,  \text{and} \,\,\, x_2 \geq \eta_{j_0+1} + \tilde{d}\}.
\end{array}\right.
\end{array}
$$
Using Lemma \ref{thm:main} on the box $\Omega_{i_0,i_0+1;j_0,j_0+1}$,
the DDM solution is found to be
$$ u = u_0 \beta_{\nearrow;\xi_{i_0+1,j_0+1}}
+ u_{\rightarrow} \beta_{\uparrow;\eta_{j_0+1}}
+ u_{\uparrow} \beta_{\rightarrow;\xi_{i_0+1}} + u_{\nearrow}.$$
\begin{remark}
When the source transfer is in the corner directions, the source for the subdomain
problem also has a specific  form, e.g., the source \eqref{eq:rhs} could be represented as
\begin{align} \label{eq:transf_form_xy}
        \Big( &\Phi_{1; i_0+1,j_0+1} (u_\uparrow;    \beta_{\rightarrow; i_0+1}) 
       + \Phi_{2; i_0+1,j_0+1} (u_\rightarrow; \beta_{\uparrow; j_0+1}) \nonumber\\
        &\Phi_{1; i_0+1,j_0+1} (u_0; \beta_{\nearrow; i_0+1}) 
       + \Phi_{2; i_0+1,j_0+1} (u_0; \beta_{\nearrow; j_0+1})
       \Big) \chi_{\nearrow; i_0+1,j_0+1},
\end{align}
where the solutions of other subdomain problems and their first order derivatives are again involved but the wave number $k$ is not.
Similar forms exist for other corner directional source transfers.  Also such source transfer in the corner directions 
could be generalized to the case of more subdomains. 
\end{remark}

With the source transfer in horizontal, vertical and corner directions as shown above,
the solution to the PML problem $\P_{\Omega}$  with the source $f_{i_0,j_0}$ could be constructed by induction, 
the details are omitted here. 
For the general source $f = \sum_{i = 1,\ldots,\Nbx \atop j = 1,\ldots,\Nby}  f_{i,j}$, the solution process for each
source $f_{i,j}$ can be handled concurrently, thus we propose
 an additive overlapping DDM for the PML problem as below.

\begin{algorithm_}[Additive overlapping DDM for the PML problem  $\P_{\Omega}$ with the source $f$] \label{alg}
$\,$
\begin{itemize}
\item Set $\{u^{0}_{i,j}\} = 0$ in $\R^2$ for $i=1,2,\ldots,N_1,$ $j =1,2,\ldots,N_2$.
\item Step 1: solve the PML problem $\P_{\Omega_{i,j}}$ with the source $f_{i,j}$
  \begin{equation} \label{eq:alg_1}
    \L_{i,j} u^{1}_{i,j} = f_{i,j}, \quad in \;\;\R^2,
  \end{equation}
for $ i=1,2,\ldots,N_1,\; j =1,2,\ldots,N_2$. 
\item 
  For Step $s = 2, 3, \ldots, N_1 + N_2$: solve
  \begin{align}  \label{eq:alg_2}
    \L_{i,j} u_{i,j}^{s} =
    &  \;\Psi_{\leftarrow; i+1,j} (u_{i+1, j}^{s-1}) + \Psi_{\rightarrow; i-1,j} (u_{i-1, j}^{s-1})   \nonumber \\
    & + \Psi_{\downarrow; i,j+1} (u_{i, j+1}^{s-1}) + \Psi_{\uparrow; i,j-1} (u_{i, j-1}^{s-1}) \nonumber \\
    & + \Psi_{\swarrow; i+1,j+1} (u_{i+1, j+1}^{s-2}) + \Psi_{\searrow; i-1,j+1} (u_{i-1, j+1}^{s-2}) \nonumber \\
    & + \Psi_{\nwarrow; i+1,j-1} (u_{i+1, j-1}^{s-2}) + \Psi_{\nearrow; i-1,j-1} (u_{i-1, j-1}^{s-2}),\quad in \;\;\R^2,
\end{align}
for $ i=1,2,\ldots,N_1,\; j =1,2,\ldots,N_2$. 
\item The DDM solution for  $\P_{\Omega}$ with the source $f$ is then given by
  \begin{equation} \label{eq:sum}
    u_{\text{DDM}} =  \sum\limits_{i = 1,\ldots,\Nbx \atop j = 1,\ldots,\Nby} 
      \beta_{0; i, j} u_{i,j} 
  \end{equation}
  with $u_{i,j} = \sum\limits_{s=1, \cdots, N_1+N_2} u^s_{i,j}.$
\end{itemize}
\end{algorithm_}

\begin{theorem}
$u_{\rm \text{DDM}}(f) = u(f)$ where $u(f)$ is the solution of the PML problem $\P_{\Omega}$ with the source $f$ and $u_{\rm \text{DDM}}(f)$ is the corresponding
DDM solution constructed by Algorithm \ref{alg}.
\end{theorem}
\emph{Proof:}
Obviously, 
$u_{\rm \text{DDM}}(f) = \sum_{i = 1,\ldots,\Nbx \atop j = 1,\ldots,\Nby} u_{\rm \text{DDM}}(f_{i,j}).$ 
Hence to show that $u_{\rm \text{DDM}}(f)$ is the solution to the problem $\P_{\Omega}$ with the source $f$,
we only need to show that for any $i_0 \in \{1,\ldots,\Nbx\}$, $j_0 \in \{1,\ldots,\Nby\}$,
 $u_{\rm \text{DDM}}(f_{i_0,j_0})$ is the solution to $\P_{\Omega}$ with the source $f(i_0,j_0)$, denoted as $u(f_{i_0,j_0})$.

\begin{figure}[ht!]	
	\begin{center} 

		\begin{minipage}{.36\textwidth}
			\centering
			\includegraphics[width=\textwidth]{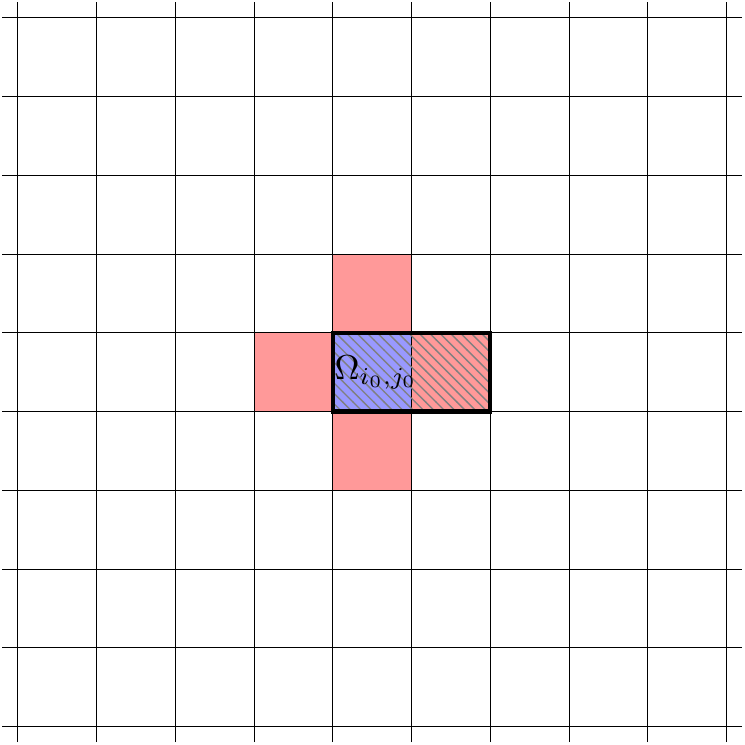}
		\end{minipage}
		\hspace{.5cm}
		\begin{minipage}{.36\textwidth}
			\centering
			\includegraphics[width=\textwidth]{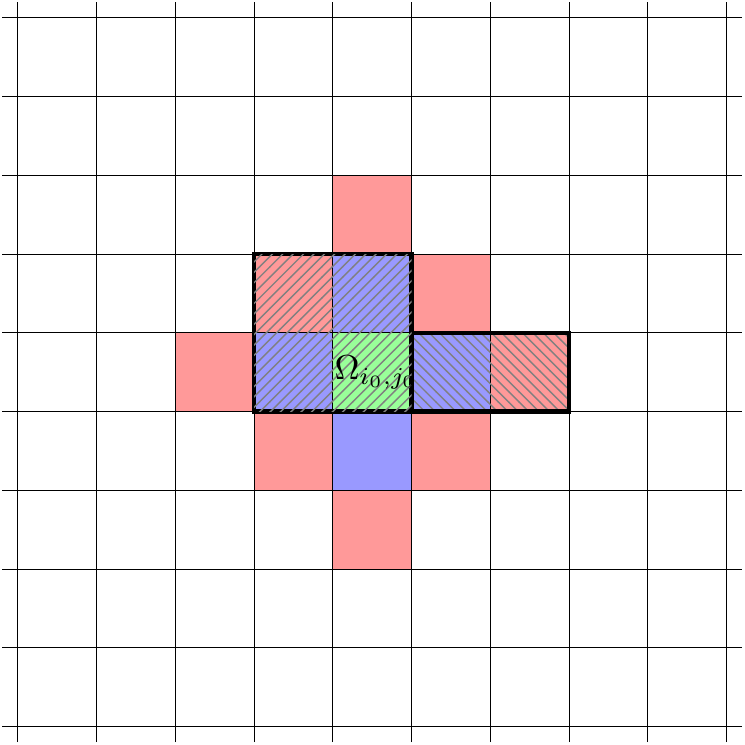}
		\end{minipage}
		\vspace{.25cm}
		\\
		\begin{minipage}{.36\textwidth}
			\centering
			Step 2
	    \end{minipage}	
		\begin{minipage}{.36\textwidth}
			\centering
			Step 3
		\end{minipage}
		\vspace{.5cm}
		\\
		\begin{minipage}{.36\textwidth}
			\centering
			\includegraphics[width=\textwidth]{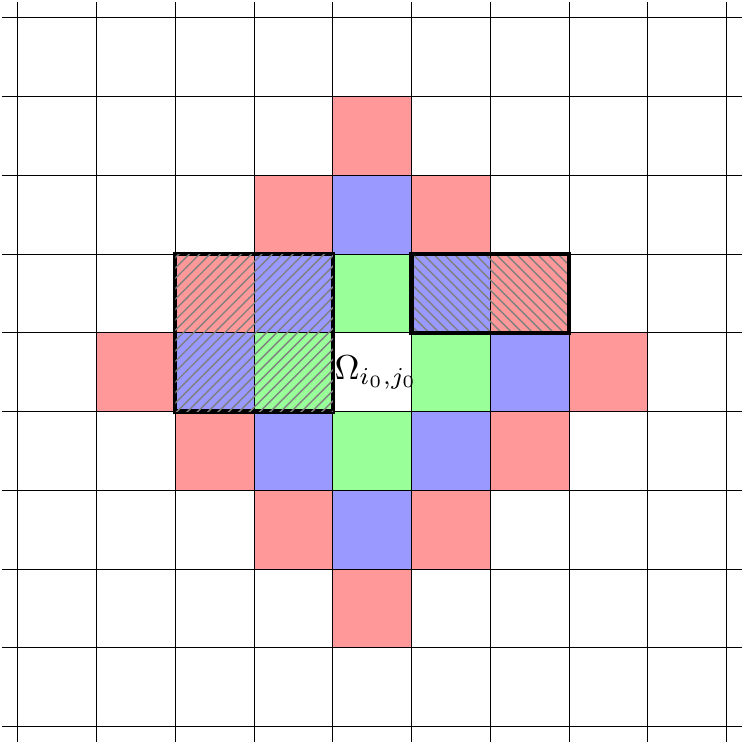}
		\end{minipage}
		\hspace{.5cm}
		\begin{minipage}{.36\textwidth}
			\centering
			\includegraphics[width=\textwidth]{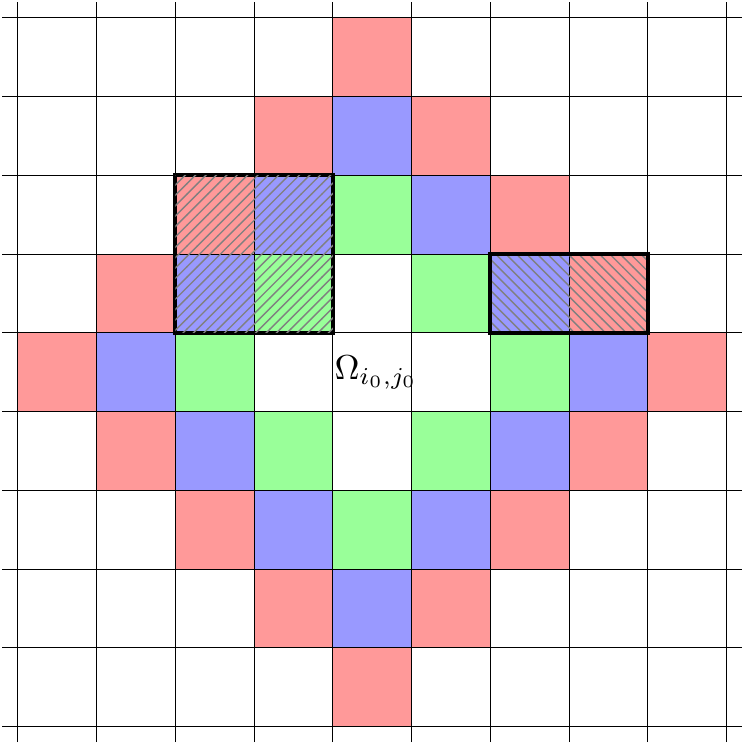}
		\end{minipage}
		\vspace{.25cm}
		\\			
		\begin{minipage}{.36\textwidth}
			\centering
			Step 4
		\end{minipage}	
		\begin{minipage}{.36\textwidth}
			\centering
			Step 5
		\end{minipage}
		\vspace{.5cm}	
		\\
		\begin{minipage}{.7\textwidth}
			\centering
			\includegraphics[width=\textwidth]{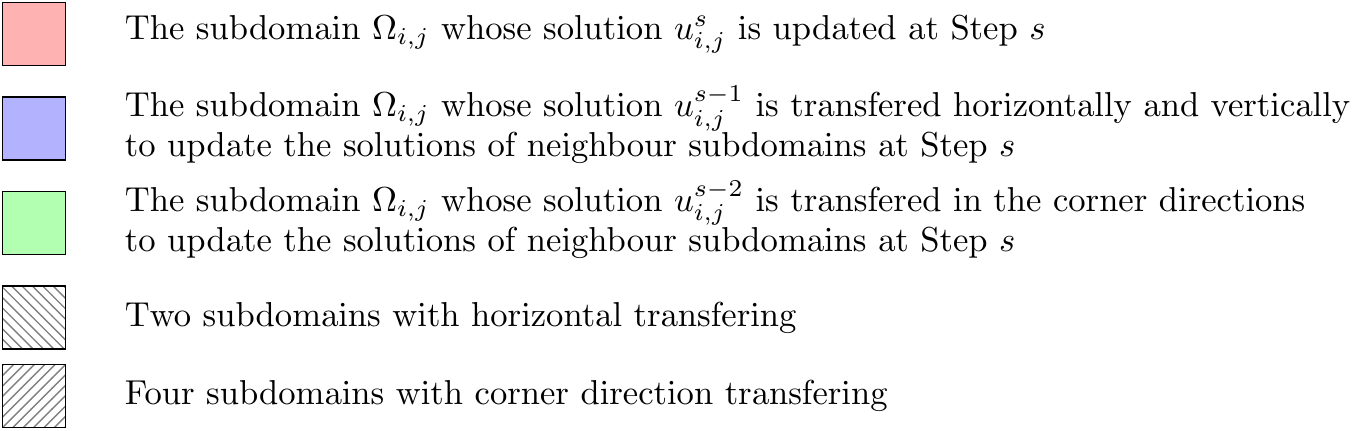}
		\end{minipage}

		\caption{The updating process of the subdomain solutions at Steps $2,3,4,5$ in the additive overlapping DDM for the PML problem $\P_{\Omega}$ with the source $f(i_0,j_0)$.}
			\label{fig:steps}
	\end{center}
\end{figure}

 At Step 1, only the solution of subdomain $\Omega_{i_0,j_0}$, $u^1_{i_0,j_0}$, need to be updated with \eqref{eq:alg_1} ($u^1_{i,j}=0$ for $i\ne i_0$ or $j\ne j_0$ due to zero source).

  At Step 2, based on \eqref{eq:alg_2}, on the subdomains $\Omega_{i,j}$ with $|i-i_0|+|j-j_0|+1 = 2$, $u_{i,j}^{2}$'s are computed with horizontal or vertical source transfers. 
  For instance, on the subdomain $\Omega_{i_0+1,j_0}$, omitting the zero terms in the RHS of \eqref{eq:alg_2}, 
  \begin{align} \label{eq:alg_horz}
    \L_{i_0+1,j_0} u_{i_0+1,j_0}^{2} 
    & =   \Psi_{\rightarrow; i_0,j_0} (u_{i_0, j_0}^{1}).
  \end{align}
  Using the analysis of horizontal and vertical  directional source transfers,
  we have that
  \begin{equation}
    \sum\limits_{r = 1,2 \atop |i-i_0|+|j-j_0|+1 \leq 2} u_{i,j}^r \beta_{0;i,j} = u(f_{i_0,j_0}),
        \qquad \text{in} \bigcup\limits_{|i-i_0|+|j-j_0|+1 \leq 2} \Omega_{i,j}
  \end{equation}
  as is shown in Figure \ref{fig:steps} (top-left).
  Note that there is no reflections opposite to the propagating direction -- let us take
  the subdomain $\Omega_{i_0+1,j_0}$ for example, with the horizontal updating \eqref{eq:alg_horz},
  $u_{i_0+1,j_0}^{2}$ is a right going ($\rightarrow$) wave solution, thus by Lemma \ref{thm:main},  $\Psi_{\Box; i_0+1,j_0}(u^2_{i_0+1,j_0}) = 0$,
  for $\Box = \leftarrow, \nwarrow, \swarrow$.

  At Step 3, 
  based on \eqref{eq:alg_2}, on  the subdomains $\Omega_{i,j}$ with $|i-i_0|+|j-j_0|+1 \leq 3$,  $u_{i,j}^{2}$'s are computed.
  However, since there is no reflections opposite to the propagating direction,
  only those on the subdomains $\Omega_{i,j}$ with $|i-i_0|+|j-j_0|+1 = 3$ really need to be updated.
  The solutions on the subdomains $\Omega_{i_0 \pm  2, j_0}$ and $\Omega_{i_0 , j_0 \pm  2}$
  are updated with  horizontal or vertical source transfers, respectively;
  the solutions on the subdomains $\Omega_{i_0 \pm  1, j_0 \pm 1}$ are updated with horizontal, vertical and corner directional source transfers.
  For instance, on the subdomain $\Omega_{i_0-1,j_0+1}$,
  omitting the zero terms in the RHS of \eqref{eq:alg_2},
  \begin{align}  \label{eq:ij0_step3}
  \L_{i_0-1,j_0+1} u_{i_0-1,j_0+1}^{3} 
  & =   \Psi_{\leftarrow; i_0,j_0+1} (u_{i_0, j_0+1}^{2})   
   + \Psi_{\uparrow; i_0-1,j_0} (u_{i_0-1, j_0}^{2})  + \Psi_{\nwarrow; i_0,j_0} (u_{i_0, j_0}^{1}).
  \end{align}
 Using the analysis of horizontal, vertical and corner directional source transfers,
  we have that
  \begin{equation}
    \sum\limits_{r = 1,\ldots,3 \atop |i-i_0|+|j-j_0|+1 \leq 3} u_{i,j}^r \beta_{0;i,j} = u(f_{i_0,j_0}),
    \qquad \text{in} \bigcup\limits_{|i-i_0|+|j-j_0|+1 \leq 3} \Omega_{i,j}
  \end{equation}
  as is shown in Figure \ref{fig:steps} (top-right).  
  Again, there is no reflections opposite to the propagating direction
  for the horizontal, vertical and corner directional transfers:
  the case of horizontal and vertical solution updating are similar to Step 2;
  for the case of corner direction  solution updating,
  let us take
  the subdomain $\Omega_{i_0-1,j_0+1}$ for example, with the corner directional  updating \eqref{eq:ij0_step3},
  $u_{i_0-1,j_0+1}^{2}$ is a up-left going ($\nwarrow$) wave solution, thus by Lemma \ref{thm:main}, $\Psi_{\Box; i_0-1,j_0+1}(u^2_{i_0-1,j_0+1}) = 0$
  for $\Box = \rightarrow, \downarrow, \nearrow, \swarrow, \searrow$.
  
  Repeat the procedure (see the bottom panel of Figure \ref{fig:steps}), we always have 
  \begin{equation}
    \sum\limits_{r = 1,\ldots,s \atop |i-i_0|+|j-j_0|+1 \leq s} u_{i,j}^r \beta_{0;i,j} = u(f_{i_0,j_0}),
     \qquad \text{in} \bigcup\limits_{|i-i_0|+|j-j_0|+1 \leq s} \Omega_{i,j}.
  \end{equation}
 Hence we finally have $ u_{\rm \text{DDM}}(f_{i_0,j_0}) = u(f_{i_0,j_0})$ after Step $N_1+N_2$.
 \qed

\section{The additive overlapping DDM in bounded truncated domains}
%
%

In this section, the additive overlapping DDM in a bounded truncated domain is proposed 
based on the DDM in the full space of $\R^2$ in the preceding section. 
Again the constant medium problem is assumed, 
and an error estimate 
of the method is established. 
First some notations are introduced as in \cite{Chen2013a}. For any bounded domain
$U \in \R^2$ with Lipschitz boundary $\Gamma$, we define the following norms:
$$||u||_{H^1(U)} = \Big(||\nabla u||^2) + ||ku||^2_{L^2(U)} \Big)^{1/2},$$
$$|||u|||_{H^1(U)} = \Big(||\nabla u||^2) + d_U^{-2}||u||^2_{L^2(U)} \Big)^{1/2},$$
$$\DD |v|^2_{\frac{1}{2}, \Gamma} = \int_{\Gamma}\int_{\Gamma} \frac{|v(x) - v(x')|^2}{|x-x'|^2} ds(x) ds(x'),$$
$$||v||_{H^{1/2}(\Gamma)} = \Big(d^{-1}_U||\nabla v||^2_{L^2(\Gamma)}) + |v|^2_{\frac{1}{2},\Gamma} \Big)^{1/2},$$
where $d_U = \text{diam}(U)$. Obviously,  
\begin{equation}
  \label{eq:H12}
  ||v||_{H^{\half}(\Gamma)} \leq
  (|\Gamma| d_U^{-1})^{\half} ||v||_{L^{\infty}(\Gamma)}
  + |\Gamma| \,\, ||\nabla v||_{L^{\infty}(\Gamma)}, \qquad \forall v \in W^{1, \infty}(\Gamma).
\end{equation}
Using the scaling argument and trace theorem we have
for any $v \in H^{\half}(\Gamma)$, 
\begin{equation} \label{eq:Hdnorm}
C_1 \frac{|U|^{\half}}{|\Gamma|} ||v||_{H^{\half}(\Gamma)} \leq \inf\limits_{ \phi|_{\Gamma} = v
\atop  \phi \in H^1(U)} |||\phi|||_{H^1(U)} \leq  C_2 \frac{|U|^{\half}}{|\Gamma|} ||v||_{H^{\half}(\Gamma)}, 
\end{equation}
where $C_1, C_2 > 0$ are constants independent of $d_U$.

\subsection{The DDM for the truncated PML problem}

Given any rectangular domain $\Omega = [-l_1,l_1]\times[-l_2,l_2]\in\R^2$, let us define an extended box $\Omega^{\PML} = [-l_1-d,l_1+d]\times[-l_2-d,l_2+d]$ with $d>\tilde d$.
Then the weak formulation of the truncated PML problem associated with $\Omega$ (denoted by $\widehat{\P}_{\Omega}$)  is defined by: for $f \in H^1(\Omega)^{\prime}$, 
find $\hat{u} \in H^1_0(\Omega^{\PML})$ such that
\begin{equation}
  \label{eq:PMLweak0}
  (A_{\Omega} \nabla \hat{u}, v) - k^2 (J_{\Omega}\hat{u}, v) = - \langle J_{\Omega}f, v \rangle, \qquad \forall\, v \in H^1(\Omega^{\PML}).
\end{equation}
Define $L = \max(l_1 + d,l_2+d)$. we consider  solve the above truncated PML problem  $\widehat{\P}_{\Omega}$ with the source $f$ using the additive overlapping DDM.

It is proved in \cite[Lemma 3.4]{Chen2013a}, that
for sufficiently large $\sigma_0 d \ge 1$, where $\sigma_0 = \max_{t \in \R} \tilde{\sigma}(t)  $,
the sesquilinear form associated with the above truncated PML problem  \eqref{eq:PMLweak0} satisfies the inf-sup condition: there exists a positive constant 
$\mu^{-1} \leq C k^{\f32}$ such that
\begin{equation}
  \label{eq:wellpose}
  \sup\limits_{\psi \in H_0^1(\Omega^{\PML})} \frac{|(A_{\Omega} \nabla \phi, \nabla \psi) - k^2 (J_{\Omega}\phi, \psi)|}{||\psi||_{H^1(\Omega^{\PML})}} 
  \ge \mu ||\phi||_{H^1(\Omega^{\PML})},
 \qquad \forall\, \phi \in H_0^1(\Omega^{\PML}).
\end{equation}
The solution  of the truncated PML problem $\widehat{\P}_{\Omega}$, denoted as  $\hat{u}$ ,
is an approximation of the solution $u$ of  the PML problem ${\P}_{\Omega}$, and the following error estimate has been proven,
\begin{equation}
\label{eq:errPML}
|| u - \hat{u} ||_{H^1(\Omega^{\PML})} \le C k^2 (1+kL)^2 e^{- \frac{1}{2} k \gamma_0 \bar{\sigma} }, 
\end{equation}
where $ \gamma_0 = \frac{d}{\sqrt{d^2 + 4 ((b_1-a_1)+(b_2-a_2) + d)^2}}$, and $\DD \bar{\sigma} = \int_{0}^{d} \tilde{\sigma}(t) dt$.
It is  clear that $\hat{u}$ is zero on $\p \Omega^{\PML}$, while ${u}$ is not. 
However, as is shown in the proof of \cite[Lemma 3.4]{Chen2013a}, 
we have the following estimate for that
\begin{align}
  |{u}(\bx)| \leq C k^{1/2} \epml ||f||_{H^1(\Omega)^{\prime}}, \qquad \forall\, \bx \in \p \Omega^{\PML}, \label{eq:estUonpB} \\
  |\nabla {u}(\bx)| \leq C k^{3/2} \epml ||f||_{H^1(\Omega)^{\prime}}, \qquad \forall\, \bx \in \p \Omega^{\PML}.
  \label{eq:estGUonpB}
\end{align}
The truncated PML problem that we solved in this 
paper is for each subdomain $\Omega_{i,j}$ ($\Omega_{i,j} \subset \Omega)$, 
therefore we use the same $\gamma_0 = \frac{d}{\sqrt{d^2 + 4(l_1+l_2 + d)^2}}$ for all truncated subdomain PML problems in the rest of the paper.

Denote $\widehat{\L}_{i,j}$ as the operator associated with the truncated subdomain PML problem $\widehat{\P}_{\Omega_{i,j}}$.
Substitute PML problem operator $\L_{i,j}$ with $\widehat{\L}_{i,j}$, 
solution $u_{\Box; i,j}$ with $\hat{u}_{\Box; i,j}$
in Algorithm \ref{alg}, 
we have the correspondingly revised DDM for the truncated PML problem in $\R^2$, 

\begin{algorithm_}[Additive overlapping DDM for the truncated PML problem $\widehat\P_{\Omega}$ with the source $f$] \label{alg2}
$\,$
\begin{itemize}
\item Set $\{\hat{u}^{0}_{i,j}\} = 0$ in  $\Omega_{i,j}^{PML}$ for $i=1,2,\ldots,N_1,$ $j =1,2,\ldots,N_2$.
\item Step 1: solve the truncated PML problem $\widehat\P_{\Omega_{i,j}}$ with the source $f_{i,j}$
  \begin{equation} \label{eq:talg_1}
    \widehat\L_{i,j} \hat{u}^{1}_{i,j} = f_{i,j}, \qquad  in \;\;\;\Omega_{i,j}^{PML},
  \end{equation}
  for $ i=1,2,\ldots,N_1,\; j =1,2,\ldots,N_2$.

\item 
  For Step $s = 2, 3, \ldots, N_1 + N_2$: solve
 \begin{align} \label{eq:talg_2}
    \widehat{\L}_{i,j} \hat{u}_{i,j}^{s} 
    = &\,   \Psi_{\leftarrow; i+1,j} (\hat{u}_{i+1, j}^{s-1}) + \Psi_{\rightarrow; i-1,j} (\hat{u}_{i-1, j}^{s-1})   \nonumber \\
    &  + \Psi_{\downarrow; i,j+1} (\hat{u}_{i, j+1}^{s-1}) + \Psi_{\uparrow; i,j-1} (\hat{u}_{i, j-1}^{s-1}) \nonumber \\
    &  + \Psi_{\swarrow; i+1,j+1} (\hat{u}_{i+1, j+1}^{s-2}) + \Psi_{\searrow; i-1,j+1} (\hat{u}_{i-1, j+1}^{s-2}) \nonumber \\
    &  + \Psi_{\nwarrow; i+1,j-1} (\hat{u}_{i+1, j-1}^{s-2}) + \Psi_{\nearrow; i-1,j-1} (\hat{u}_{i-1, j-1}^{s-2}), \quad  in \;\;\;\Omega_{i,j}^{PML},  
\end{align}
for $ i=1,2,\ldots,N_1,\; j =1,2,\ldots,N_2$. 

\item The DDM solution for  $\widehat{P}_{\Omega}$ with the source $f$ is then given by

  \begin{equation} \label{eq:tsum}
    \hat{u}_{\text{DDM}}(f) =  \sum\limits_{i = 1,\ldots,\Nbx \atop j = 1,\ldots,\Nby} 
      \beta_{0; i, j} \hat{u}_{i,j} 
  \end{equation}
  with $\hat{u}_{i,j} = \sum\limits_{s=1, \cdots, N_1+N_2} \hat{u}^s_{i,j}.$
\end{itemize}
\end{algorithm_}

\begin{remark}
The additive overlapping DDM algorithms (Algorithms \ref{alg} and \ref{alg2}) can be straightforwardly generalized to solve the three-dimensional PML and truncated PML problems,  where there are now a total of $3^3-1=26$ directions for source transfer and $N_1+N_2+N_3$ iteration steps ($N_3$ denotes the number of partitions in the $z$-direction).
\end{remark}

\subsection{Error estimate}
Let us first only consider the case of only the source $f_{i_0, j_0} \neq 0$. 
Then it holds that the subdomain solution $u^s_{i,j}$ is non-zero only at Step $s = |i-i_0|+|j-j_0|+1$,
while the subdomain solution $\hat{u}^s_{i,j}$ is non-zero at Steps $ s \geq |i-i_0|+|j-j_0|+1$.

\begin{figure}[ht!]	
	\begin{center}
		\includegraphics[width=.45\textwidth]{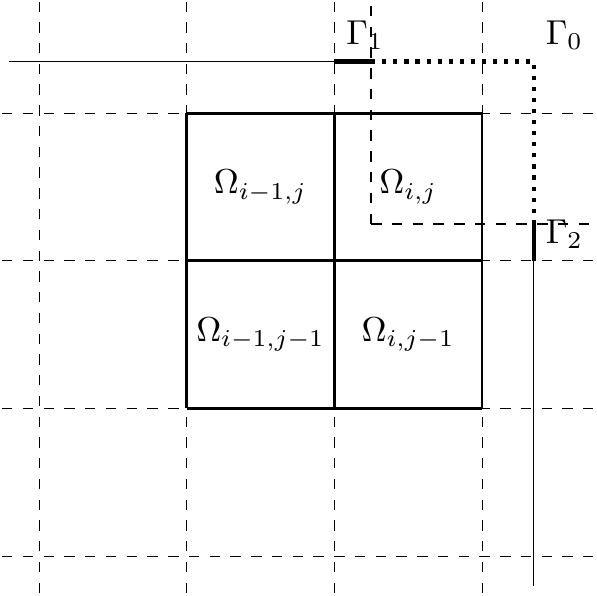}
		\caption{$\Gamma_0$, $\Gamma_1$, $\Gamma_2$.}	
		\label{fig:uOnPML}
	\end{center}
\end{figure}

%
%

\begin{lemma} \label{lemma:UonPML}
  Assume that only the source $f_{i_0, j_0} \neq 0$ and $\sigma_0 d$ be sufficiently large. For $i=1,\dots,N_1$, $j=1,\dots,N_2$, 
  it holds
  \begin{equation}
    \label{eq:UonPML}
    ||{u}_{i,j}^s||_{H^{1/2}(\p \Omega_{i,j}^{\PML})}
    \leq C k^{\half}(1+kL) \epml ||f_{i_0,j_0}||_{H^1(B_l)^{\prime}},           
  \end{equation}
  for $s = |i-i_0|+|j-j_0|+1$.
\end{lemma}

\emph{Proof:}
  First we prove that 
  \begin{align} 
     |u_{i,j}^{s}(\bx)| \leq& C k^{\half} \epml ||f_{i_0,j_0}||_{H^1(B_l)^{\prime}},  &\forall\, \bx \in \p \Omega_{i,j}^{\PML}, \label{eq:maxUOnPML} \\
   |\nabla u_{i,j}^{s}(\bx)| \leq& C k^{\f32} \epml ||f_{i_0,j_0}||_{H^1(B_l)^{\prime}},  &\forall\, \bx \in \p \Omega_{i,j}^{\PML}. \label{eq:maxGuOnPML} 
  \end{align}
  Consider the case $i \ge i_0 + 1$, $j \ge j_0 + 1$ for example.
  Denote by $u_{i_0,i; j_0,j}$  the solution to the PML problem $\P_{\Omega_{i_0,i; j_0,j}}$ with $f_{i_0,j_0}$ as the source, then
  we have  by \eqref{eq:estUonpB}-\eqref{eq:estGUonpB}
  \begin{align} 
    |u_{i_0,i; j_0,j}| \le C k^{\half} \epml ||f_{i_0,j_0}||_{H^1(B_l)^{\prime}}, \qquad\forall\,\bx \in \p B_{i_0,i; j_0,j}.\label{eq:box_i0ij0j}\\
    |\nabla u_{i_0,i; j_0,j}| \le C k^{\f32} \epml ||f_{i_0,j_0}||_{H^1(B_l)^{\prime}}, \qquad\forall\,\bx \in \p B_{i_0,i; j_0,j}.\label{eq:box_i0ij0j_g}
  \end{align}

  Define
  $\Gamma_0 = \p \Omega_{i,j}^{\PML} \cap \{x_1 \ge \xi_i + d, x_2 \ge \eta_j + d\}$,
  $\Gamma_1 = \{ \xi_i \le x_1 \le \xi_i + d, x_2 = \eta_{j+1} + d \}$,
  $\Gamma_2 = \{ x_1 = \xi_{i+1} + d, \eta_j \le x_2 \le \eta_{j+1} + d\}$, as shown in Figure \ref{fig:uOnPML}.
  By Lemma \ref{thm:main}, $u_{ i,j}^{s}$ is zero on $\p \Omega_{i_0,i; j_0,j}$ except for $\Gamma_0 \cup \Gamma_1 \cup \Gamma_2$, moreover,
  by \eqref{eq:corner_u},
  $$
  \begin{array}{rll}
    u_{ i,j}^s =& u_{i_0,i; j_0,j},& \qquad{\rm on} \;\;\Gamma_0,\\
    u_{ i,j}^s =& u_{i_0,i; j_0,j} (1-\beta_{\rightarrow; i_0}) ,&\qquad{\rm on}\;\; \Gamma_1,\\
    u_{ i,j}^s =& u_{i_0,i; j_0,j} (1-\beta_{\uparrow;  j_0}),&\qquad  {\rm on}\;\; \Gamma_2,
  \end{array}
  $$
which imply  \eqref{eq:maxUOnPML}-\eqref{eq:maxGuOnPML} together with  \eqref{eq:box_i0ij0j}-\eqref{eq:box_i0ij0j_g}.

  By \eqref{eq:H12} we then obtain
  \begin{equation}
    ||{u}_{i,j}^s||_{H^{1/2}(\p \Omega_{i,j}^{\PML})}
    \leq C \max\limits_{x \in \p \Omega_{i,j}^{\PML}} (|{u}_{i,j}^s| + L |\nabla {u}_{i,j}^s|),
  \end{equation}
  which further gives us \eqref{eq:UonPML} by using \eqref{eq:maxUOnPML} and \eqref{eq:maxGuOnPML}.
\qed

%
%
\begin{lemma} \label{thm:err0}
	 Assume that only the source $f_{i_0, j_0} \neq 0$ and $\sigma_0 d$ be sufficiently large. For $i=1,\dots,N_1$, $j=1,\dots,N_2$, 
	it holds
	\begin{equation} \label{eq:trun_err}
	||u_{ i,j}^s - \hat{u}_{ i,j}^s||_{H^{1}(\Omega_{i,j}^{\PML})} 
           \leq C k^{\f32 s + \half}  (1+kL)^2 \epml ||f_{i_0,j_0}||_{H^1(\Omega)^{\prime}}.
	\end{equation}
\end{lemma}
\emph{Proof:}
  Let the source of the PML problem $\P_{\Omega_{i,j}}$ at Step $s$ in \eqref{eq:alg_1}-\eqref{eq:alg_2}
  of the DDM Algorithm \ref{alg}
  be denoted as ${f}_{i,j}^s$, 
  and the source for the truncated PML problem $\widehat{\P}_{\Omega_{i,j}}$ at Step $s$ in 
  \eqref{eq:talg_1}-\eqref{eq:talg_2}
  of the DDM algorithm \ref{alg2} be denoted as $\hat{f}_{i,j}^s$.  

  First we will prove the following error estimate for the source in the DDM Algorithm by induction
  	\begin{equation} \label{eq:trun_f_err}
	||f_{ i,j}^s - \hat{f}_{ i,j}^s||_{H^{1}(\Omega_{i,j}^{\PML})} 
	\leq C k^{\f32 s - 1}  (1+kL)^2 \epml ||f_{i_0,j_0}||_{H^1(\Omega)^{\prime}}.
	\end{equation}  
    Obviously \eqref{eq:trun_f_err} holds for $s=1$ since 
    $f_{ i,j}^1 = \hat{f}_{ i,j}^1$, suppose it holds for $s=2,\ldots,t-1$, we will prove it 
    also holds for $s=t$. 
  By the analysis of source transfers in horizontal, vertical and corner directions, 
  we know that the sources of the subdomain problems have other forms that involve
 the solutions of other subdomain problems and their first-order derivatives but not $k$.
  This indicates
  \begin{align}
  || {f}_{ i,j}^t - \hat{f}_{ i,j}^t ||_{H^{-1}(\Omega_{i,j}^{\PML})} &\leq 
     \sum\limits_{(i',j') = (i\pm 1,j), (i,j \pm 1)}
    || {u}_{ i',j'}^{t-1} - \hat{u}_{ i',j'}^{t-1} ||_{H^{1}(\Omega_{i',j'}^{\PML})} \nonumber \\
    & \quad + \sum\limits_{(i',j') = (i\pm 1,j+1),(i\pm 1,j-1)}
    || {u}_{ i',j'}^{t-2} - \hat{u}_{ i',j'}^{t-2} ||_{H^{1}(\Omega_{i',j'}^{\PML})}.
  \label{eq:trun_f_on_u}      
  \end{align}
  
  Take $ || {u}_{i+1,j}^{t-1} - \hat{u}_{i+1,j}^{t-1} ||_{H^{1}(\Omega_{i+1,j}^{\PML})} $ for example, 
  by \eqref{eq:Hdnorm}, there exists a lifting function $\phi \in H^1(\Omega_{i+1,j}^{\PML})$ that $\phi = {u}_{i+1,j}^{t-1}$
  on $\p \Omega_{i+1,j}^{\PML}$, and $|||\phi|||_{H^1(\Omega_{i+1,j}^{\PML})} \leq C || {u}_{i+1,j}^{t-1} ||_{H^{\half}(\p \Omega_{i+1,j}^{\PML})} $. 
  Note that $ {u}_{i+1,j}^{t-1} - \hat{u}_{i+1,j}^{t-1} - \phi \in H^1_0(\Omega_{i+1,j}^{\PML}) $, then by \eqref{eq:wellpose}, 
  we have
  \begin{align*}
|| {u}_{i+1,j}^{t-1} - \hat{u}_{i+1,j}^{t-1} - \phi ||_{H^{1}(\Omega_{i+1,j}^{\PML})} 
&\leq \mu^{-1} \bigg( \sup\limits_{\psi \in H^1_0(\Omega_{i+1,j}^{\PML}) } \frac{|(A_{\Omega_{i+1,j}}\nabla \phi, \nabla \psi) - k^2(J_{\Omega_{i+1,j}}\phi, \psi )|}{||\psi||_{H^1(\Omega_{i+1,j}^{\PML})}} \\
&\quad + || {f}_{i+1,j}^{t-1} - \hat{f}_{i+1,j}^{t-1} ||_{H^{-1}(\Omega_{i+1,j}^{\PML})}    \bigg)\\
&\leq C \mu^{-1} \Big( ||\phi||_{H^1(\Omega_{i+1,j}^{\PML})}
  + || {f}_{i+1,j}^{t-1} - \hat{f}_{i+1,j}^{t-1} ||_{H^{-1}(\Omega_{i+1,j}^{\PML})} \Big),
\end{align*}
therefore we have
  \begin{align*}
    || {u}_{i+1,j}^{t-1} - \hat{u}_{i+1,j}^{t-1} ||_{H^{1}(\Omega_{i+1,j}^{\PML})} 
    &\leq C k^{\f32}  || {f}_{i+1,j}^{t-1} - \hat{f}_{i+1,j}^{t-1} ||_{H^{-1}(\Omega_{i+1,j}^{\PML})} 
     + Ck^{\f32}(1+kL) ||| \phi |||_{H^{1}(\Omega_{i+1,j}^{\PML})} \\
    &\leq C k^{\f32}  || {f}_{i+1,j}^{t-1} - \hat{f}_{i+1,j}^{t-1} ||_{H^{-1}(\Omega_{i+1,j}^{\PML})} \\
    & \quad + C k^{\f32} (1+kL) || {u}_{i+1,j}^{t-1} ||_{H^{\half}(\p \Omega_{i+1,j}^{\PML})}.
  \end{align*}
  Since ${u}_{i,j}^{s} =  0$ for $|i-i_0|+|j-j_0| + 1 \neq s$, and we have \eqref{eq:UonPML} for $|i-i_0|+|j-j_0| + 1 = s$  in Lemma \ref{lemma:UonPML}, thus
   \begin{align*}
   || {u}_{i+1,j}^{t-1} - \hat{u}_{i+1,j}^{t-1} ||_{H^{1}(\Omega_{i+1,j}^{\PML})} 
   &\leq C k^{\f32}  || {f}_{i+1,j}^{t-1} - \hat{f}_{i+1,j}^{t-1} ||_{H^{-1}(\Omega_{i+1,j}^{\PML})} \\
   &\quad + C k^{2} (1+kL)^2 \epml || f_{i_0,j_0} ||_{H^1(\Omega)^{\prime}}.
   \end{align*} 
  Similar estimates also  hold for the other terms in the right hand side of \eqref{eq:trun_f_on_u}.
  Therefore we have
  \begin{align*}
    || {f}_{ i,j}^{t} - \hat{f}_{ i,j}^{t} ||_{H^{-1}(\Omega_{i,j}^{\PML})} &\leq 
	 \sum\limits_{(i',j') = (i\pm 1,j), (i,j \pm 1)}
	C k^{3/2} || {f}_{ i',j'}^{t-1} - \hat{f}_{ i',j'}^{t-1} ||_{H^{1}(\Omega_{i',j'}^{\PML})} \nonumber \\
	& \quad + \sum\limits_{(i',j') = (i\pm 1,j+1),(i\pm 1,j-1)}
	C k^{3/2} || {f}_{ i',j'}^{t-2} - \hat{f}_{ i',j'}^{t-2} ||_{H^{1}(\Omega_{i',j'}^{\PML})} \nonumber \\
    &\quad +   C k^{2} (1+kL)^2 \epml || f_{i_0,j_0} ||_{H^1(\Omega)^{\prime}} \\
    &\leq    C k^{\f32 t - 1}  (1+kL)^2 \epml ||f_{i_0,j_0}||_{H^1(\Omega)^{\prime}},
  \end{align*}
  thus we have \eqref{eq:trun_f_err} by induction.

  Next we will prove \eqref{eq:trun_err}. For the case $s=1$, \eqref{eq:trun_err} holds by \eqref{eq:errPML},  and for the case $s>1$,  it also holds by using standard argument again 
  with  \eqref{eq:UonPML} and  \eqref{eq:trun_f_err} .
\qed

In the above analysis, we deal with the case of only  $f_{i_0, j_0} \neq 0$. 
For the general case, $f = \sum\limits_{i, j} f_{i, j}$, 
$u(f) = \sum\limits_{i, j} u(f_{i, j})$, 
thus we have the following result by Lemma \ref{thm:err0}.
\begin{theorem} \label{thm:err}
	Assume  $\sigma_0 d$ is sufficiently large, then it holds
	\begin{equation} \label{eq:errDDM}
	||u - \hat{u}_{\rm \text{DDM}}||_{H^{1}(\Omega^{\PML})} \leq C k^{\f32 N_1 + \f32 N_2 + \half} (1+kL)^2 \epml ||f||_{H^1(\Omega)^{\prime}},
	\end{equation}
        where $u$ is the solution to the PML problem $\P_{\Omega}$ with the source $f$, 
        and $\hat{u}_{\rm \text{DDM}}$ is the DDM solution for the corresponding truncated PML problem constructed by Algorithm \ref{alg2}.
\end{theorem}

\section{Numerical experiments}

We will  numerically test our additive overlapping DDM in bounded truncated domains using Algorithm \ref{alg2}.
First, the convergence of the discrete DDM solutions is tested.
From Theorem \ref{thm:err}, the truncated DDM solution in the continuous case is a good 
approximation to the solution of the original problem  \eqref{eq:helm} when appropriate PML medium
parameters are chosen. 
With the subdomain problems being discretized by finite difference or
finite element method, it is expected that the discretization error usually dominate the total 
error, and to verify this, we test our method  
for two and three dimensional problems with constant medium. 
Second, our DDM will  also be tested as an iterative  solver (i.e., Algorithm  \ref{alg2} could be used as an iterative method in which the total number of iterations/steps is determined by certain convergence/stopping criterion) or as a preconditioner for global discrete systems for
large wave number problems.  We test Algorithm  \ref{alg2} for  both some constant medium and 
layered medium problems to compare their performance.

Five points finite difference scheme in two dimensions and seven points finite difference scheme in three dimensions are respectively
taken to discretize the problems (the Laplacian operator) on uniform rectangular meshes,  both of them are obviously second order accurate. 
We have implemented a parallel version of DDM algorithm \ref{alg2} using Message Passing Interface (MPI) and the local 
subdomain problems   are solved with the direct solver ``MUMPS'' \cite{MUMPS}. 
The supercomputer ``Tianhe-2 Cluster'' located in Tianjin, China is used for our numerical tests, each node of 
the cluster has two 2.2GHz Xeon E5-2692 processors with 12 cores, and 
64G memory. 

\subsection{Convergence tests of the discrete DDM solutions}
In this subsection, the proposed DDM is tested for the problems with  constant medium  --
with the wave number and the number of subdomains both being fixed, 
the mesh resolution is increased to test the convergence of
the discrete solutions.

\subsubsection{2D example}
In this example, we solve a two-dimensional Helmholtz problem with a constant wave number
$k/2\pi = 25$. 
The computational domain is $B_L=[-L,L]\times [-L,L]$ with $L = 1/2$, and the interior domain without PML is
$B_l=[-l.l]\times [l,-l]$ with $l = 25/56$.
Denote $q$ as the mesh density which is the number of nodes per wave length, 
a serie of meshes are used in the computation where the mesh density is
approximately  $q \approx 22,45,90,180,270$, respectively.
The domain $B_L$ is partitioned into $5 \times 5$ subdomains,
and thus  the total number of steps used for constructing the DDM solution is $5+5=10$.
%
The source is chosen as  
\begin{equation*}
f = \frac{16 k^2}{\pi^3}  e^{- (\frac{4 k}{\pi})^2 ((x-x_0)^2+(y-y_0)^2) },
\end{equation*}
where $(x_0, y_0) = (0.09, 0.268)$, thus it is almost supported by four subdomains
$\Omega_{3,4}$, $\Omega_{3,5}$, $\Omega_{4,4}$ and $\Omega_{4,5}$.

\begin{figure}[ht!]	
	\begin{center}
		\includegraphics[width=.6\textwidth]{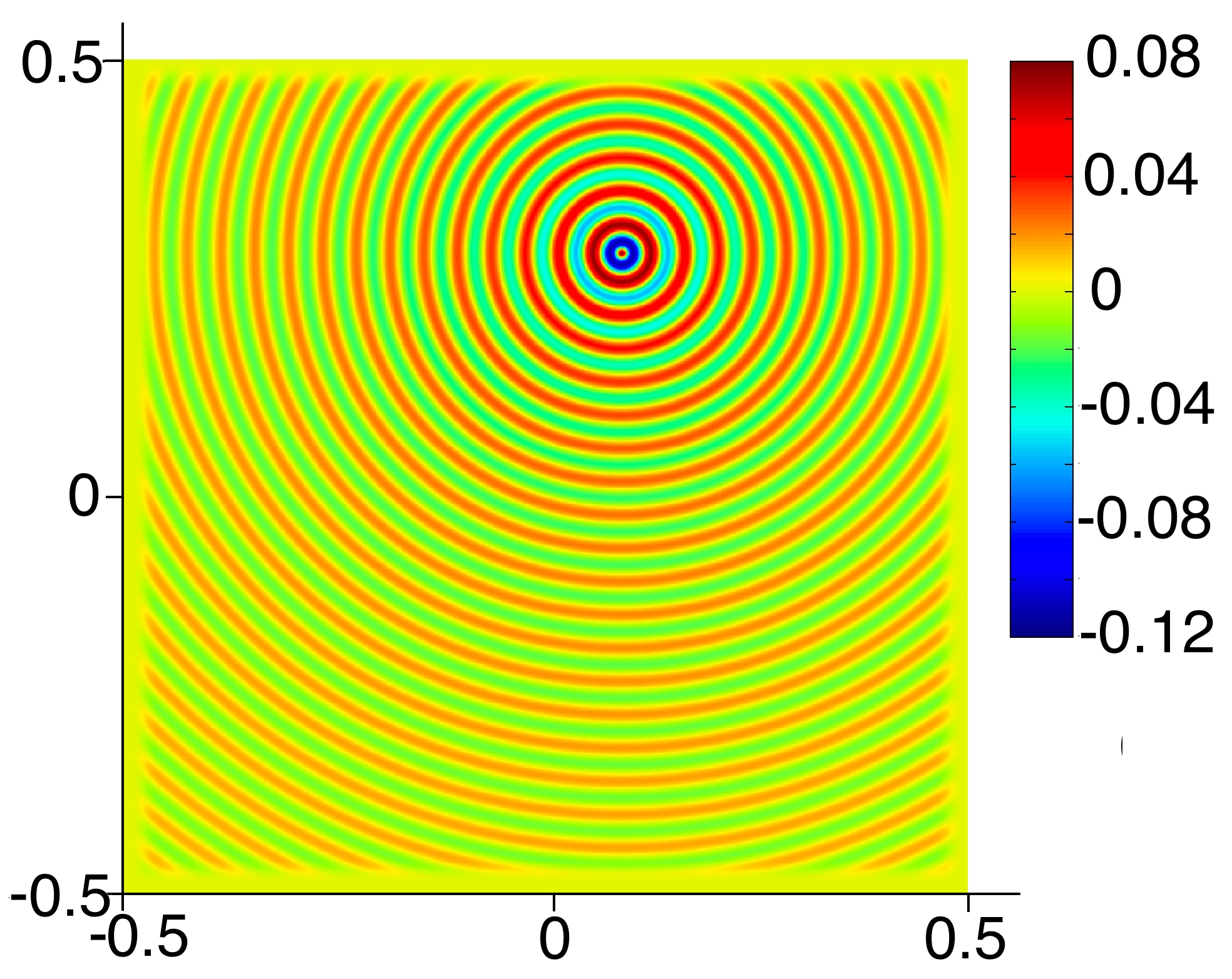}
		\caption{The real part of the simulated DDM solution on the mesh of size 2240$^2$ for the 2D constant medium problem.}	
		\label{fig:num_sol1}
	\end{center}
\end{figure}

Table \ref{tab:num2D}
reports  the errors and convergence rates of the discrete DDM solutions  for this 2D example, and Figure \ref{fig:num_sol1} shows the real part of the simulated 
DDM solution on the mesh of 6720$^2$.
It is observed that the errors of the truncated DDM solutions \eqref{eq:errDDM} are very small with appropriate 
medium parameters for PML and the finite difference discretization errors do dominate the total errors.
We obtain as expected the optimal second order convergence of the errors measured by the $L^2$ norm along the refinement of the meshes. 
Note that we also have super-convergence of the  $H^1$ error.
Such behavior was also observed in numerical tests for Poisson problem \cite{Ng2009} and  proved on rectangular domains \cite{Strikwerda2004}.

\begin{table}[ht!]
	\centering
	\begin{tabular}{|r|c|r|c|r|c|}
		\hline
		Mesh  & Local Size & $L^2$ Error  & Conv. & $H^1$  Error & Conv. \\
		Size & without PML &                     &   Rate    &                      &  Rate     \\
		\hline\hline
		560$^2$   & 100$^2$   & 3.13$\times 10^{-3}$ &       & 5.47$\times 10^{-1}$ &       \\
		1120$^2$  & 200$^2$   & 7.68$\times 10^{-4}$ & 2.0    & 1.36$\times 10^{-1}$ & 2.0    \\
		2240$^2$  & 400$^2$   & 1.95$\times 10^{-4}$ & 2.0    & 3.40$\times 10^{-2}$ & 2.0    \\
		4480$^2$  & 800$^2$   & 4.85$\times 10^{-5}$ & 2.0    & 8.50$\times 10^{-3}$ & 2.0    \\
		6720$^2$  & 1200$^2$  & 2.16$\times 10^{-5}$ & 2.0    & 3.79$\times 10^{-3}$ & 2.0     \\
		\hline
	\end{tabular}
	\caption{The errors and convergence rates of the DDM solutions for the 2D constant medium problem.} \label{tab:num2D}
\end{table}

\subsubsection{3D example}
Next we solve a three-dimensional Helmholtz problem with a constant wave number
$k/2\pi = 10$. 
The computational domain is still $B_L$ with $L = 1/2$, but the interior domain without PML is 
$B_l$ with $l = 3 / 8$.
A series of meshes is used in the computation where the mesh density is
approximately  $q \approx 8, 9.5, 13, 16$, respectively.
The domain $B_L$ is partitioned into $3 \times 3 \times 3$ subdomains,
and the total number of steps used for constructing the DDM solutions is $3+3+3=9$.
%
The source is chosen as  
\begin{equation*}
f = \frac{64 k^3}{\pi^{9/2}}  e^{- (\frac{4 k}{\pi})^2 ((x-x_0)^2+(y-y_0)^2+(z-z0)^2 ) },
\end{equation*}
where $(x_0, y_0, z_0) = (0.12, 0.133, 0.125)$, thus it is almost supported by eight subdomains
$\Omega_{i,j,m}$ where $i=2,3$, $j=2,3$, $m=2,3$.

\begin{figure}[ht!]	
	\begin{center}
		\includegraphics[width=.65\textwidth]{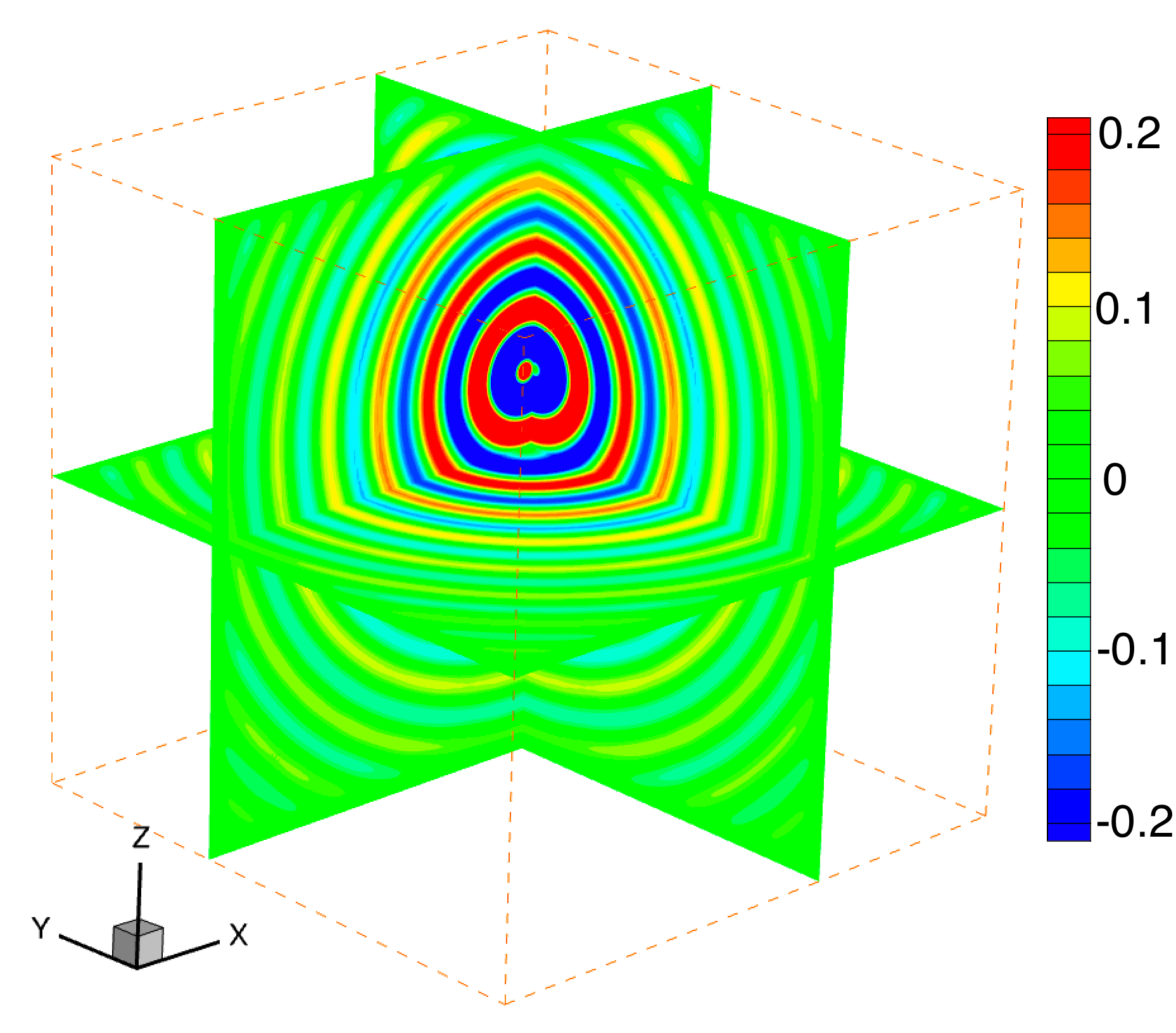}
		\caption{The real part of the simulated DDM solution on the mesh of size 160$^3$ for the 3D constant medium problem.}	
		\label{fig:num_sol2}
	\end{center}
\end{figure}

Table \ref{tab:num3D}
reports  the errors and convergence rates of the discrete DDM solutions for this 3D example, and Figure \ref{fig:num_sol2} shows the real part of the simulated 
solution on the mesh of 160$^3$.
Again, the finite difference errors dominate the total errors, 
and we obtain the optimal second order convergence for both $L^2$ and $H^1$ errors as expected.

\begin{table}[ht!]
	\centering
	\begin{tabular}{|r|c|r|c|lc|}
		\hline
		Mesh  & Local Size & $L^2$ Errors               & Conv. & $H^1$ Errors              & Conv. \\
		Size& without PML &                     &   Rate    &                      &  Rate     \\
		\hline     \hline                                                                                          
		80$^3$    & 20$^3$                      & 2.66$\times 10^{-2}$ &       & 1.82$\times 10^{0}$ &       \\
		96$^3$    & 24$^3$                      & 1.81$\times 10^{-2}$ & 2.1   & 1.25$\times 10^{0}$ & 2.1   \\
		128$^3$   & 32$^3$                      & 9.99$\times 10^{-3}$ & 2.1   & 6.94$\times 10^{-1}$ & 2.0   \\
		160$^3$   & 40$^3$                      & 6.39$\times 10^{-3}$ & 2.0   & 4.45$\times 10^{-1}$ & 2.0   \\
		\hline
	\end{tabular}
	\caption{The errors and convergence rates of the DDM solutions for the 3D constant medium problem.} \label{tab:num3D}
\end{table}

\subsection{Performance tests with the DDM as an iterative solver or a preconditioner  for the global discrete system}

The DDM Algorithm \ref{alg2} designed for the constant medium Helmholtz problem
also could be used as an iterative solution method for both constant and variable medium problems.
More iterations/steps than $N_1+N_2$ is then usually needed for the iterative solver
to reach certain relative error tolerances. We  first test the DDM method as an iterative solver 
for some 2D constant medium and layered medium problems and then  as a preconditioner for the global discrete system and compare their performance.

A single shot in the subdomain $\Omega_{1,1}$ is taken as the source
with the position being $(x^s, y^s)$, where 
\begin{equation}
x^s = x_0 + \Delta \lx/4,\qquad y^s = y_0 + \Delta \ly/3.
\end{equation}
Let $h_x$ and $h_y$
 denote the mesh sizes in $x$ and $y$ directions, respectively.
 The shape of the shot is an approximate $\delta$ function, defined by
 \begin{equation}\DD f_{i,j} = \frac{1}{h_x h_y}\delta(x^s - i h_x) \delta (y^s - j h_y ).\end{equation}

The size of the subdomain problems is fixed in all of the following tests, each non-overlapping subregion is
of size $300 \times 300$, and the PML layer is 30 points width at each side of the subdomain.
The number of subdomains is increased from $2\times2$ to $32 \times 32$.
Denote by $\nIterDDM$ the total number of  iterations.  
The number of DDM solve defined by
\begin{equation}
\nSolveDDM =\DD \frac{\nIterDDM}{N_1+N_2}
\end{equation}
is used to measure the effectiveness of the iterative solver.
It is always highly desired that the number of  used iterations $\nIterDDM$  for convergence is proportional to $N_1+N_2$
as the number of subdomains grows,
in that case, $\nSolveDDM$ remains almost the same.
%
The tolerance for relative residuals is set to be 10$^{-8}$ as the stopping criterion for all tests.
Note that we also let each subdomain use exactly one core in the following tests.

\subsubsection{As an iterative solver}
The Algorithm \ref{alg2} is firstly tested as an iterative solver 
for the global discrete system.

First we test the additive overlapping DDM solver for a constant medium problem on the 
square domain $[0,1]\times[-1,0]$, with different frequencies. The results on the numbers of used iterations  and the running times are reported in Table \ref{tab:iter_const}.
As we can see, the number of DDM solves $\nSolveDDM$ grows as the number of the subdomains grows from 
$2\times2$ to $8\times8$, but remains around 1.92 when the number of the subdomains grows further.

\begin{table}[ht!]
	\centering
	\begin{tabular}{|r|c|r|r|r|r|r|}
		\hline
		Mesh & $N_1 \times N_2$ & Frequency                            & $\nIterDDM$ & $\nSolveDDM$ & \multicolumn{1}{c|}{Total}& Time/Iter \\
		\multicolumn{1}{|c|}{Size}&                    & \multicolumn{1}{c|}{$\omega / 2\pi$} &            &                  & Time (sec)           &    \multicolumn{1}{c|}{ (sec)} \\
		\hline  
		\hline  
		600$^2$    & 2 $\times$ 2       & 55                                 & 5          & 1.25              & 2.0      &    0.40            \\ 
		1200$^2$  & 4 $\times$ 4       & 105                                 & 13        & 1.63              & 4.7      &     0.36             \\   
		2400$^2$  & 8 $\times$ 8       & 205                                & 29         & 1.81              & 11.1       &    0.38            \\  
		4800$^2$  & 16 $\times$ 16     & 405                                & 61         & 1.91              & 24.4       &   0.40             \\  
		9600$^2$  & 32 $\times$ 32     & 805                                & 123        & 1.92              & 55.9        &  0.45             \\    
		19200$^2$  & 64 $\times$ 64    & 1605                                & 247        & 1.93              & 117.2      &   0.47              \\    
		\hline
	\end{tabular}
	\caption{The performance of the DDM as an iterative solver  for the 2D constant medium problem
		with the subdomain problem size being fixed.} \label{tab:iter_const}
\end{table}

\begin{figure}[ht!]
\centerline{
\includegraphics[width=.52\textwidth]{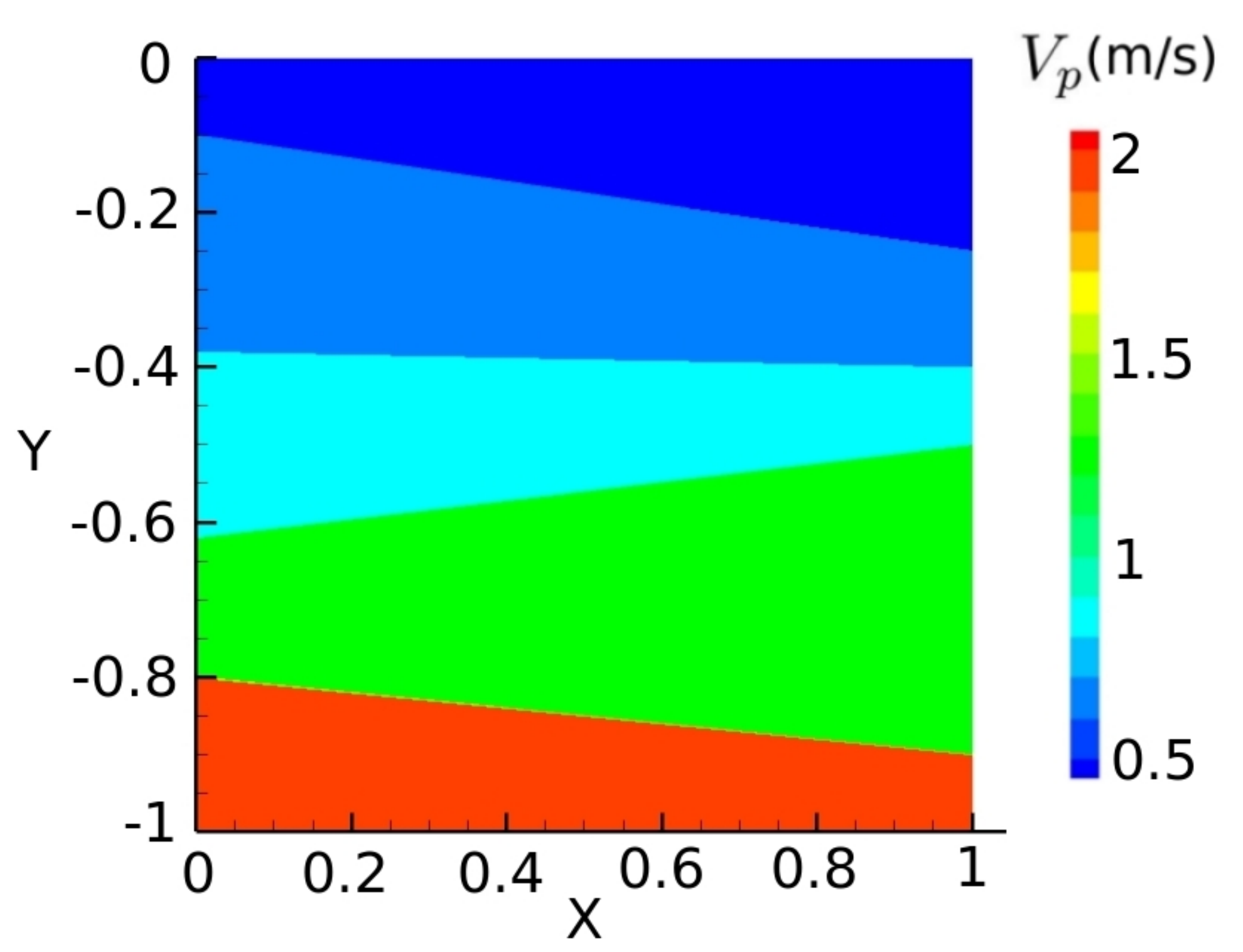} 
\includegraphics[width=.51\textwidth]{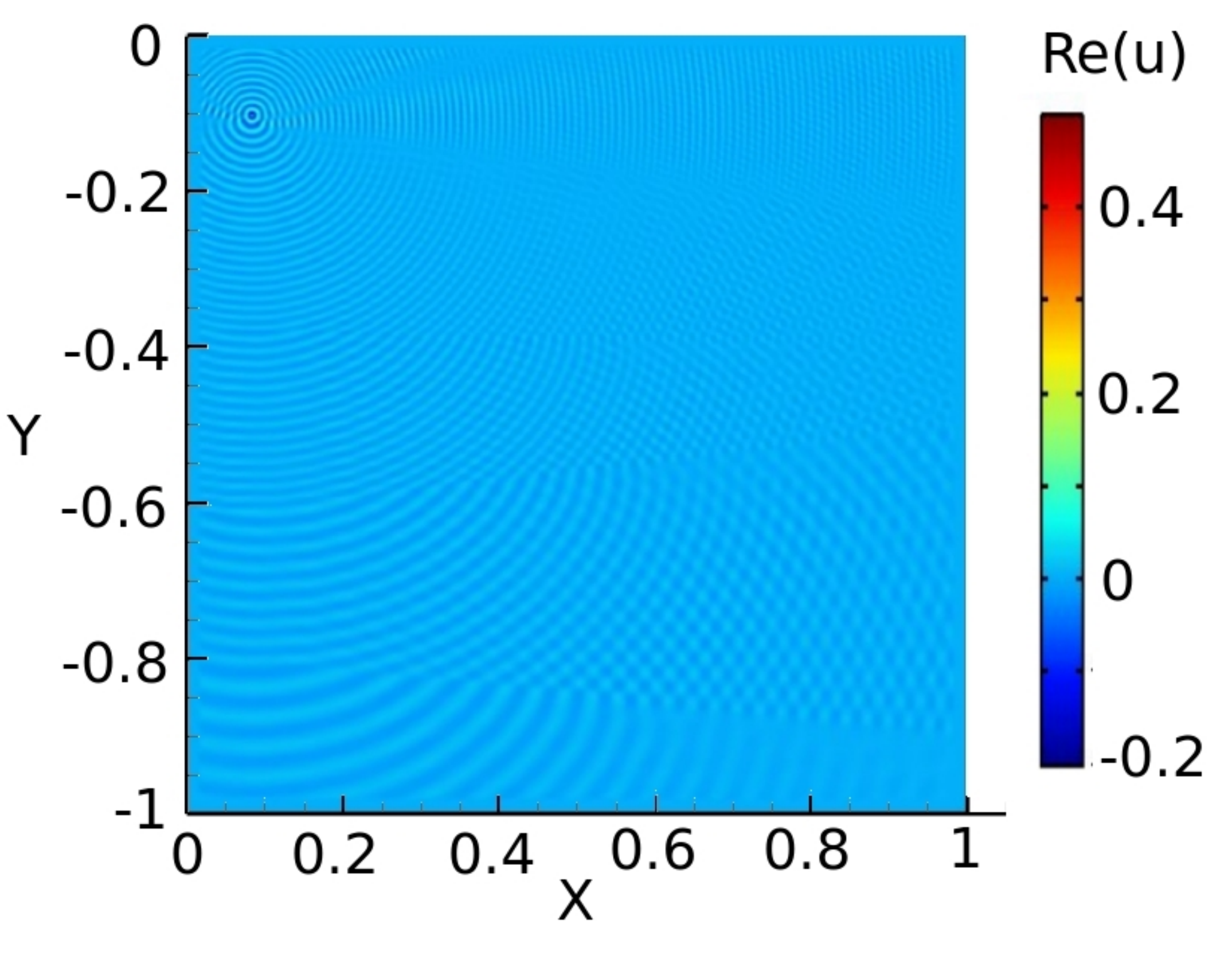} 
}
\caption{Left: the 2D layered medium velocity profile; 
	right: the real part of discrete solution of the problem on the mesh of  size $2400^2$. 
}
\label{fig:5layer}
\end{figure}

%
%
%
%
Next we test the DDM solver for a layered medium problem on the 
square domain $[0,1]\times[-1,0]$, as shown in Figure \ref{fig:5layer} (left) for the velocity profile.
For the purpose of illustration,  
the approximate solution of the problem produced by using the uniform mesh of size $1600^2$ is shown in Figure \ref{fig:5layer} (right).
The results on the  numbers of used iterations  and the running  times  are reported in Table \ref{tab:iter_layer}.
As we can see,  the number of DDM solves almost does not change  as the number of the subdomains grows,
and it requires about 2.3 times iterations to converge compared with the constant medium problem. 

\begin{table}[ht!]
\centering
\begin{tabular}{|r|c|r|r|r|r|r|}
  \hline
	Mesh & $N_1 \times N_2$ & Frequency                            & $\nIterDDM$ & $\nSolveDDM$ & \multicolumn{1}{c|}{Total}& Time/Iter \\
		\multicolumn{1}{|c|}{Size}&                    & \multicolumn{1}{c|}{$\omega / 2\pi$} &            &                  & Time (sec)           &    \multicolumn{1}{c|}{ (sec)} \\
\hline  
\hline  
 600$^2$    & 2 $\times$ 2       &   27.5                                   & 16         &  4.00             & 5.7           &    0.36           \\  
 1200$^2$  & 4 $\times$ 4       &   52.5                                   & 38         &  4.75            & 13.2           &   0.35          \\ 
 2400$^2$  & 8 $\times$ 8       &   102.5                                  & 73         &  4.56             & 27.8            &    0.38       \\
 4800$^2$  & 16 $\times$ 16     &   202.5                                  & 146        &  4.56             & 62.4           & 0.43           \\
 9600$^2$  & 32 $\times$ 32     &   402.5                                  & 291        &  4.55             & 133.1           &  0.46           \\
  19200$^2$  & 64 $\times$ 64     &   802.5                                  & 582        &  4.55             & 278.8          &   0.48           \\   
\hline
\end{tabular}
\caption{The performance of the DDM as an iterative solver  for the 2D layered medium problem
with the subdomain problem size being fixed.} \label{tab:iter_layer}
\end{table}

Based on the results on the  time cost per DDM iteration with up to $64\times 64 =4096$ cores (equivalent to the number of subdomains) reported in Tables \ref{tab:iter_const} and \ref{tab:iter_layer}, it is  easy to see that the proposed additive overlapping DDM is indeed highly scalable  and achieve excellent 
parallel efficiency.

\subsubsection{As a preconditioner}
Next we use Algorithm \ref{alg2} as a preconditioner
to the flexible GMRES (FGMRES) solver for the global discrete system. The number of subdomains is set to $16\times 16$ as an example.
Note that when using Algorithm \ref{alg2} as a preconditioner, the number of DDM steps $K$ in each
preconditioning solve is to be chosen, we test different $K$'s for both constant medium  and
layered medium cases, and the results are reported in Tables \ref{tab:K-all} and \ref{tab:K-all-layered}.
Since the major computational cost during the iterations of solving the linear system is in the back substitute of LU factorization of each
subdomain, so the number of back substitutes, denoted as $\nLocalSolve$, is also
presented  in the tables.

In the constant medium case, both the total running time  and $\nLocalSolve$ basically decrease   
as $K$ increases from 1 to $N_1+N_2=16+16=32$.  It is also easy to see that $\nLocalSolve$ for all cases of $K$ are 
 larger than $\nIterDDM$ ($= 61$) used by the DDM iterative solver with $16\times 16$ subdomains (see Table \ref{tab:iter_const}),
 and so does the total time cost (24.4 seconds for the DDM iterative solver). The DDM preconditioner with the largest $K=32$ achieves
almost similar performance  as the DDM iterative solver.

In the layered  medium case, $\nLocalSolve$ only slightly increase (one exception)
 from 148 to 160 as $K$ increases.
 When $K=1$, although $\nLocalSolve$ is relatively small, the extra cost of GMRES takes more time and the total running time is the largest.
 When $K=2$, the best  $\nLocalSolve$ ($=146$) is reached and in fact it equals exactly, $\nIterDDM$, the number of iterations needed by the DDM iterative solver with $16\times 16$ subdomains (see Table \ref{tab:iter_layer}), but the total running time is larger (71.3 vs. 62.4 seconds).
 For $K \ge 3$,  all the time costs are close to 62.4 seconds, which indicates that the performance differences between the  DDM preconditioner 
 and the DDM iterative solver are small.

 To summarize, 
 it is observed by our experiments that when the additive overlapping DDM is used as a preconditioner, large $K$ should be taken, and 
 its performance is slightly worse than to be used as an iterative solver.

\begin{table}[ht!]
\centering
\begin{tabular}{|c|r|r|r|}
  \hline
    $K$   & $\nIterGMERES$ & $\nLocalSolve$ & Total Time (sec)\\
  \hline\hline
 1     & 91             & 91             & 44.0 \\
 2     & 41             & 82             & 39.8 \\
                      3     & 27             & 81             & 37.8 \\
   5     & 16             & 80             & 35.0 \\
 10    & 8              & 80             & 33.3 \\
 15    & 5              & 75             & 30.7 \\
                      20    & 4              & 80             & 32.9 \\
                     32    & 2              & 64             & 26.6 \\
  \hline
\end{tabular}
\caption{The performance of the  DDM  as a preconditioner to the global FGMRES solver for the 2D constant medium problem.
The number of subdomains is $16\times 16$.}
\label{tab:K-all}
\end{table}

\begin{table}[ht!]
\centering
\begin{tabular}{|c|r|r|r|}
  \hline
         $K$   & $\nIterGMERES$ & $\nLocalSolve$ & Total Time (sec) \\
  \hline\hline
    1  & 148            & 148            & 83.7 \\
                      2  & 73             & 146            & 71.3 \\
                      3  & 50             & 150            & 68.5 \\
 5  & 30             & 150            & 66.1 \\
 10 & 15             & 150            & 62.3 \\
 15 & 10             & 150            & 60.8 \\
   20 & 8              & 160            & 66.9 \\
   32 & 5              & 160            & 63.1 \\
  \hline
\end{tabular}
\caption{The performance of the  DDM  as a preconditioner to the global FGMRES solver for the 2D layered medium problem.
The number of subdomains is $16\times 16$.}\label{tab:K-all-layered}
\end{table}

\section{Conclusions}

In this paper, we propose and analyze an additive overlapping DDM for solving the high-frequency Helmholtz equation.
The source transfer technique is used in horizontal, vertical and corner directions, to make the wave 
propagates  away from the source subdomain by subdomain. For the constant medium case, 
the DDM method is shown to  exactly solve the  PML problem defined in $\mathbb{R}^2$, and 
an error estimation is established for the method used for solving the PML problem in bounded truncated domains.
The DDM method could be used as an iterative solver or a preconditioner for the global discrete system in both constant or
 variable medium problems,  and as an iterative solver is more preferable.
The method is  very suitable for large-scale parallel computing as demonstrated by numerical experiments.
Although the method has been extended to three dimensions,  
 error estimation of the method for the truncated PML problem in bounded three-dimensional domains is still open
and will be studied in the next step.
In addition, we would like to remark that this DDM method could be generalized to irregular domains, with special domain decompositions  
that require each subdomain to be convex polygons.  In such case, the source transfer 
would be along the outward normal directions of each subdomain boundary. The convergence analysis and performance of 
such generalized method is also worthy of further investigation.



\end{document}